\documentclass{amsart} 
\usepackage{amssymb,latexsym} 
 
\numberwithin{equation}{section} 
 
\newtheorem{theorem}{Theorem}[section]

\newtheorem{corollary}[theorem]{Corollary} 
\newtheorem{remark}[theorem]{Remark} 
\newtheorem{lemma}[theorem]{Lemma} 
\newtheorem{example}[theorem]{Example}

\def\proof{\smallskip\noindent {\bf Proof. }} 
\def\endproof{\hfill$\square$\medskip}

\newcommand{\rem}[1]{\left\langle#1\right\rangle} 
 
\newcommand{\overunder}[2]{ 
\!\begin{array}{c} 
\scriptstyle{#1}\\[-.1in] 
-\!\!\!-\!\!\!-\\[-.1in] 
\scriptstyle{#2} 
\end{array} 
\! 
} 
\newcommand{\wideoverunder}[2]{ 
\!\begin{array}{c} 
\scriptstyle{#1}\\[-.1in] 
-\!\!\!-\!\!\!-\!\!\!-\!\!\!-\!\!\!-\!\!\!-\\[-.1in] 
\scriptstyle{#2} 
\end{array} 
\! 
} 
\newcommand{\overunderarrow}[2]{ 
\!\begin{array}{c} 
\scriptstyle{#1}\\[-.1in] 
-\!\!\!-\!\!\!\to\\[-.1in] 
\scriptstyle{#2} 
\end{array} 
\! 
} 
\newcommand{\wideoverunderarrow}[2]{ 
\!\begin{array}{c} 
\scriptstyle{#1}\\[-.1in] 
-\!\!\!-\!\!\!-\!\!\!-\!\!\!-\!\!\!-\!\!\!\to\\[-.1in] 
\scriptstyle{#2} 
\end{array} 
\! 
}

\def\AAA{\mathbb{A}}

\def\TT{\mathbb{T}} 
\def\ZZ{\mathbb{Z}}

\def\LL{\mathcal{L}}

\begin{document} 
 
\title[The Laurent phenomenon]{The Laurent phenomenon 
%\\ for exchange recurrences 
}

% Information for first author~/www/Combin/winter01/ 
\author{Sergey Fomin} 
\address{Department of Mathematics, University of Michigan, 
Ann Arbor, MI 48109, USA} \email{fomin@umich.edu}

\author{Andrei Zelevinsky} 
\address{\noindent Department of Mathematics, Northeastern University, 
 Boston, MA 02115} 
\email{andrei@neu.edu} 
 
\date{April 25, 2001} 
 
\thanks{The authors were supported in part 
by NSF grants \#DMS-0049063, \#DMS-0070685 (S.F.), and \#DMS-9971362 
(A.Z.). 
}

% General info 
\subjclass{Primary 14E05, % Rational and birational maps 
Secondary 
05E99, % None of the above, but in this section 
11B83. % Special sequences and polynomials 
%12E05% Polynomials (irreducibility, etc.) 
} 
 
 
\begin{abstract} 
A composition of birational maps given by Laurent polynomials 
need not be given by Laurent polynomials; however, sometimes---quite 
unexpected\-ly---it does.  
We suggest a unified treatment of this phenomenon,  
which covers a large class of applications.  
In particular, we settle in the affirmative a conjecture 
of D.~Gale and R.~Robinson %\cite{gale-intelligencer, guy-unsolved} 
on integrality of generalized Somos sequences, and prove the Laurent property for  
several multidimensional recurrences, confirming  
conjectures by J.~Propp, N.~Elkies, and M.~Kleber.  
\end{abstract} 
 
%{\ } 
 
%\vspace{-.9in} 
  
\maketitle 
 
\tableofcontents

\section{Introduction} 
\label{sec:intro} 
 
In this paper, we suggest a unified explanation for a number of  
instances in which certain recursively defined rational functions 
prove, unexpectedly, to be Laurent polynomials.  
We begin by presenting several instances %providing several manifestations 
of this \emph{Laurent phenomenon} established in the paper.

%%% This happens in various seemingly unrelated contexts,  
%%% each time making for a surprising discovery.  
%%% In this paper, we suggest a general treatment of most known 
%%% manifestations of the Laurent phenomenon, providing a unified explanation  
%%% for a miscellany of applications.  
% In particular, we settle in the affirmative a conjecture 
% of D.~Gale and R.~Robinson \cite{gale-intelligencer, guy-unsolved} 
% on integrality of generalized Somos sequences, and prove the Laurent property for  
% several multidimensional recurrences, confirming  
% conjectures made by J.~Propp, N.~Elkies, and M.~Kleber.  
%%% Here are three characteristic examples.  
 
\begin{example}  
\label{example:kleber} 
(The cube recurrence) 
{\rm  
Consider a 3-dimensional array 
\[ 
(y_{ijk}\,:\, (i,j,k) \in \mathcal{H})
%\ZZ,\, i+j+k\geq 0) 
\] 
whose elements satisfy the recurrence  
\begin{equation} 
\label{eq:cube-frac} 
y_{i,j,k} 
=\frac{ 
 \alpha y_{i-1,j,k}y_{i,j-1,k-1} 
+\beta y_{i,j-1,k}y_{i-1,j,k-1} 
+\gamma y_{i,j,k-1}y_{i-1,j-1,k}}{y_{i-1,j-1,k-1} 
}. 
\end{equation}
%More precisely, let 
Here $\mathcal{H}$ can be any non-empty subset of $\ZZ^3$ satisfying 
the following conditions: 
\begin{eqnarray} 
\label{eq:cube-h1} 
&&\text{if $(i,j,k)\in\mathcal{H}$, then $(i',j',k')\in\mathcal{H}$ 
whenever $i\leq i',j\leq j',k\leq k'$;} 
\\ 
\label{eq:cube-h2} 
&&\text{for any $(i',j',k')\in\mathcal{H}$, the set 
$\{(i,j,k)\in\mathcal{H}:i\leq i',j\leq j',k\leq k'\}$ is finite.}
\end{eqnarray}  
} 
\end{example} 
 
\begin{theorem} 
\label{th:kleber} 
Let
%\[ 
$H_{\rm init} = \{ (a,b,c)\in\mathcal{H}\,:\,(a-1,b-1,c-1)\notin\mathcal{H} \}$. 
%\]
For every $(i,j,k)\in\mathcal{H}$, the entry $y_{i,j,k}$ is a Laurent
polynomial with coefficients in $\ZZ[\alpha,\beta,\gamma]$ 
in the initial entries $y_{a,b,c}$, for $(a,b,c) \in H_{\rm init}$. 
%$a+b+c\in\{0,1,2\}$.  
\end{theorem} 
 
The cube recurrence (with $\alpha=\beta=\gamma=1$) 
was introduced by James Propp \cite{propp}, who was also the one to 
conjecture Laurentness in the case when 
$\mathcal{H} \subset \ZZ^3$ is given by the condition $i+j+k\geq 0$;
in this case $H_{\rm init}$ consists of all $(a,b,c) \in \mathcal{H}$
such that $a+b+c\in\{0,1,2\}$.
Another natural choice of $\mathcal{H}$ was suggested by Michael Kleber:
$\mathcal{H} = \ZZ_{\geq 0}^3$, in which case 
$H_{\rm init} = \{(a,b,c) \in \ZZ_{\geq 0}^3 : abc = 0\}$.

\begin{example} 
\label{example:gale-robinson} 
(The Gale-Robinson sequence) 
{\rm 
Let $p$, $q$, and $r$ be distinct positive integers, 
let $n=p+q+r$, and let 
the sequence $y_0,y_1,\dots$ satisfy the recurrence 
\begin{equation} 
\label{eq:gale-robinson} 
y_{k+n}=\frac{\alpha y_{k+p}y_{k+n-p}+\beta y_{k+q}y_{k+n-q} 
+\gamma y_{k+r}y_{k+n-r}}{y_{k}}\,. 
\end{equation}  
David Gale % 
and Raphael Robinson %\cite{robinson} 
conjectured (see~\cite{gale-intelligencer} and
\cite[E15]{guy-unsolved}) that every term of such a sequence  
is an integer provided $y_0=\cdots=y_{n-1}=1$ 
and $\alpha,\beta,\gamma$ are positive integers. 
Using Theorem~\ref{th:kleber},  
we prove the following stronger statement. 
} 
\end{example} 
 
\begin{theorem} 
\label{th:gale-robinson} 
As a function of the initial terms $y_0$,\dots,~$y_{n-1}$, 
every term of the Gale-Robinson sequence %{\rm(\ref{eq:gale-robinson})} 
is a Laurent polynomial with coefficients in $\ZZ[\alpha,\beta,\gamma]$.  
\end{theorem} 

We note that the special case $\alpha = \beta = \gamma =1$, $p=1$, 
$q=2$, $r=3$, $n=6$  
(resp., $r=4$, $n=7$) of the recurrence (\ref{eq:gale-robinson}) 
is the Somos-6 (resp., Somos-7) recurrence \cite{gale-intelligencer}. 
 
\begin{example} 
\label{example:octahedral}  
(Octahedron recurrence)  
{\rm  
Consider the 3-dimensional recurrence 
\begin{equation}
\label{eq:oct}
y_{i,j,k}  
=\frac{\alpha y_{i+1,j,k-1} y_{i-1,j,k-1}  
+\beta y_{i,j+1,k-1} y_{i,j-1,k-1} }{y_{i,j,k-2}} 
\end{equation} 
for an array $(y_{ijk})_{(i,j,k)\in\mathcal{H}}$ whose indexing set 
$\mathcal{H}$ is contained in the lattice 
\begin{equation}
\label{eq:lattice-mod2}
L=\{(i,j,k)\in\ZZ^3 \,:\, i+j+k\equiv 0\bmod 2\} 
\end{equation}
and satisfies the following analogues of conditions 
(\ref{eq:cube-h1})--(\ref{eq:cube-h2}): 
\begin{eqnarray} 
\label{eq:oct-h1} 
&&\text{if $(i,j,k)\in\mathcal{H}$, then $(i',j',k')\in\mathcal{H}$ 
whenever $|i'-i|+|j'-j|\leq k'-k$;} 
\\ 
\label{eq:oct-h2} 
&&\text{for any $(i',j',k')\in\mathcal{H}$, 
the set $\{(i,j,k)\in\mathcal{H}: |i'-i|+|j'-j|\leq k'-k\}$}\\
\nonumber
&&\text{is finite.} 
\end{eqnarray} 
} 
\end{example} 
 
\begin{theorem} 
\label{th:oct} 
Let $H_{\rm init} = \{ (a,b,c)\in\mathcal{H}\,:(a,b,c-2)\notin\mathcal{H} \}$.
For every $(i,j,k)\in\mathcal{H}$, the entry $y_{i,j,k}$ is a Laurent
polynomial with coefficients in $\ZZ[\alpha,\beta]$ 
in the initial entries $y_{a,b,c}$, for $(a,b,c)\in H_{\rm init}$.  
\end{theorem} 
 
The octahedron recurrence on the half-lattice 
%\[
\begin{equation}
\label{eq:half-lattice-mod2}
\mathcal{H}=\{(i,j,k)\in L \,:\,k\geq 0\}
\end{equation}
%, \, i+j+k\equiv 0\bmod 2\}
%\]  
was studied by W.~H.~Mills, D.~P.~Robbins, and  
H.~Rumsey in their pioneering work~\cite{mills-robbins-rumsey}  
on the Alternating Sign Matrix Conjecture 
(cf.\ \cite{bressoud-propp} and \cite[Section~10]{propp} 
for further references); in particular, they proved the special case 
of Theorem~\ref{th:oct} for this choice of~$\mathcal{H}$. 
%\[  
%H_{\rm init}=\{(i,j,k)\in\mathcal{H}\,:\,k\in\{0,1\}\} \ . 
%\]  

\begin{example} 
\label{example:two-term-gale-robinson} 
(Two-term version of the Gale-Robinson sequence) 
{\rm 
Let $p$, $q$, and $n$ be positive integers 
such that $p<q\leq n/2$,  
and let the sequence $y_0,y_1,\dots$ satisfy the recurrence 
\begin{equation} 
\label{eq:gale-robinson-two-term} 
y_{k+n}=\frac{\alpha y_{k+p}y_{k+n-p}+\beta y_{k+q}y_{k+n-q} 
}{y_{k}}\,. 
\end{equation}
Using Theorem~\ref{th:oct}, one can prove that this sequence also exhibits
the Laurent phenomenon.
}
\end{example} 

\begin{theorem} 
\label{th:two-term-gale-robinson} 
As a function of the initial terms $y_0$,\dots,~$y_{n-1}$, 
every term $y_m$ %of the Gale-Robinson sequence %{\rm(\ref{eq:gale-robinson})} 
is a Laurent polynomial with coefficients in $\ZZ[\alpha,\beta]$.  
\end{theorem} 

We note that in the special case  
$\alpha = \beta =1$, $p=1$, $q=2$, $n=5$ (resp., $n=4$), 
(\ref{eq:gale-robinson-two-term}) becomes 
the Somos-5 (resp., Somos-4) recurrence \cite{gale-intelligencer}.

The last example of the Laurent phenomenon 
presented in this section is of a somewhat different 
kind; it is inspired by \cite{conway-guy}. 
  
\begin{example} %(The Conway recurrence) 
\label{example:conway-scott} 
{\rm  
Let $n \geq 3$ be an integer, 
and consider a quadratic form 
$$P(x_1, \dots, x_n) = x_1^2 + \cdots + x_n^2 + 
\displaystyle \sum_{i < j} \alpha_{ij} x_i x_j \ .$$ 
Define the rational transformations $F_1,\dots,F_n$ by 
% in the variables $x_1,x_2,x_3$  
\begin{equation} 
F_i : (x_1,\dots,x_n) \mapsto  
(x_1,\dots,x_{i-1},\frac{P\bigl|_{x_i = 0}}{x_i},x_{i+1},\dots,x_n). 
\end{equation} 
} 
\end{example} 
 
\begin{theorem} 
\label{th:conway-scott} 
For any sequence of indices $i_1,\dots,i_m$, 
%where $1\leq i_k\leq n$ for every~$k$, 
the composition map 
$G=F_{i_1}\circ\cdots\circ F_{i_m}$ is given by 
\[ 
G : x=(x_1,\dots,x_n) \mapsto (G_1(x),\dots,G_n(x)), 
\] 
where $G_1,\dots,G_n$ are Laurent polynomials %in $x_1,\dots,x_n$ 
with coefficients in~$\ZZ[\alpha_{ij}:i<j]$. 
\end{theorem} 
 
%This example can be vastly generalized by introducing 
%coefficients, changing the exponents (with proper care), etc. 

\medskip
  
This paper is an outgrowth of \cite{fz-clust1}, where we 
initiated the study of a new class of commutative algebras, 
called cluster algebras, 
and established the Laurent phenomenon in that context. 
Here we prove the theorems stated above, along with a number of 
related results, 
using an approach inspired by \cite{fz-clust1}. 
The first step is to reformulate the problem in terms of 
generalized \emph{exchange patterns} (cf.\ \cite[Definition~2.1]{fz-clust1}),  
which consist of \emph{clusters} and \emph{exchanges} among them. 
The clusters are distinguished finite sets of variables,  
each of the same cardinality~$n$.  
An exchange operation on a cluster ${\bf x}$ replaces 
a variable $x \in {\bf x}$ by a new variable 
$x' = \frac{P}{x}$, where $P$ is a polynomial in the $n-1$ variables 
${\bf x}-\{x\}$.  
Each of the above theorems 
can be restated as saying that any member of the cluster obtained from an  
initial cluster ${\bf x}_0$ by a particular sequence of exchanges  
is a Laurent polynomial in the variables from~${\bf x}_0$.  
Theorem~\ref{th:conway-scott} is explicitly stated in this way; 
in the rest of examples above, the rephrasing is less straightforward. 
  
Our main technical tool is ``The Caterpillar Lemma''  
(Theorem~\ref{th:laurent-gcd}), which establishes the Laurent  
phenomenon for a particular class of exchange patterns 
(see Figure~\ref{fig:caterpillar}). 
This is a %significant 
modification of the namesake statement \cite[Theorem~3.2]{fz-clust1}, 
and its proof closely follows the argument in~\cite{fz-clust1}. 
(We note that none of the two statements is a formal consequence of another.) 
  
In most applications, including Theorems~\ref{th:kleber} 
and \ref{th:oct} above, the ``caterpillar'' patterns to which 
Theorem~\ref{th:laurent-gcd} applies, are not manifestly present within the  
original setup.  
Thus, we first complete it by creating  
additional clusters and exchanges, and then apply the Caterpillar Lemma. 
 
% The work presented here grew out of our attempts to establish 
% the Laurent phenomenon for \emph{cluster algebras,} a new class of 
% algebraic objects we introduced \cite{fz-clust1} in our study of total 
% positivity and dual canonical bases in semisimple Lie groups.  
% We then realized that the Laurent phenomenon holds in a very  
% general context, and can be proved by essentially elementary means.  
 
The paper is organized as follows.  
The Caterpillar Lemma is proved in Section~\ref{sec:caterpillar}.  
Subsequent sections contain its applications.  
In particular, Theorems~\ref{th:kleber}, \ref{th:gale-robinson}, 
\ref{th:oct}, and \ref{th:two-term-gale-robinson} are  
proved in Section~\ref{sec:lattice-based},
while Theorem~\ref{th:conway-scott} is proved
in Section~\ref{sec:homogeneous-exchange-patterns}.  
Other instances of the Laurent phenomenon treated in this paper
include generalizations of each of the following: 
Somos-4 sequences (Example~\ref{example:somos4}),
Elkies's ``knight recurrence'' (Example~\ref{example:elkies}), 
frieze patterns (Example~\ref{example:frieze-patterns}) 
and number walls (Example~\ref{example:number-walls}).  
 
We conjecture that in all instances of the Laurent phenomenon 
established in this paper, the Laurent polynomials in question  
have \emph{nonnegative} integer coefficients.  
In other contexts, similar nonnegativity conjectures were made earlier in  
\cite{fz-jams, fz-intelligencer, fz-clust1}.
 
\medskip 
 
\vbox{\noindent \textsc{Acknowledgments.}  
We thank Jim Propp for introducing us to a number of 
beautiful examples of the Laurent phenomenon, 
and for very helpful comments on the first draft of the paper. 
In particular, it was he who showed us how to deduce 
Theorem~\ref{th:two-term-gale-robinson} from Theorem~\ref{th:oct}. 
 
This paper was completed during our 
stay at the Isaac Newton Institute 
for Mathematical Sciences (Cambridge,~UK), 
whose support and hospitality are gratefully acknowledged. 
} 
 
%\section{Generalized exchange patterns} 
%\label{sec:exchange} 
 
\section{The Caterpillar Lemma} 
\label{sec:caterpillar} 
 
Let us fix an integer $n\geq 2$, and let $T$ be a tree whose edges are  
labeled by the elements of the set $[n]=\{1,2,\dots,n\}$, 
so that the edges emanating from each vertex receive 
different labels. 
By a common abuse of notation, we will sometimes denote by $T$ 
the set of the graph's vertices. 
We 
%say that two vertices $t,t'\in T$ are 
%\emph{$k$-adjacent,} and 
will write $t \overunder{k}{} t'$ if vertices 
$t,t'\in T$ are joined by an edge labeled by~$k$. 
 
 From now on, let $\AAA$ be a unique factorization domain (the ring of integers 
$\mathbb{Z}$ or a suitable polynomial ring would suffice for most 
applications).  
Assume that a nonzero polynomial $P\in\AAA[x_1,\dots,x_n]$, 
not depending on~$x_k\,$, is associated  
with every edge $t \overunder{k}{} t'$ in~$T$.  
We will write $~t \overunder{}{P} t'~$ or $~t \overunder{k}{P} 
t'~$, and call $P$ the \emph{exchange polynomial} associated with the 
given edge. 
The entire collection of these polynomials is called a 
\emph{generalized exchange pattern} on~$T$.  
(In \cite{fz-clust1}, we introduced a much narrower notion of an  
\emph{exchange pattern;} hence the terminology.)  
 
We fix a root vertex $t_0\in T$, and introduce the \emph{initial  
cluster} ${\bf x} (t_0)$ of $n$ independent variables $x_1(t_0),\dots,x_n(t_0)$. 
To each vertex $t\in T$, we then associate a \emph{cluster} 
${\bf x} (t)$ consisting of $n$  
elements $x_1(t),\dots,x_n(t)$ of the field %$\FFcal$ 
of rational functions $\AAA(x_1(t_0),\dots,x_n(t_0))$. 
The elements $x_i(t)$ are uniquely determined by the following 
\emph{exchange relations,} for 
every edge $t \overunder{k}{P} t'$: 
\begin{eqnarray} 
\label{eq:exchange1} 
&&\text{$x_i(t)=x_i(t')$\quad for any $i\neq k$;}\\%[.1in] 
\label{eq:exchange2} 
&&\text{$x_k(t)\,x_k(t')=P({\bf x}(t))$.} 
\end{eqnarray} 
(One can recursively compute the $x_i(t)$'s, moving away from the  
root. Since the exchange polynomial $P$ does not depend on $x_k$,  
the exchange relation (\ref{eq:exchange2}) does not change if we apply  
it in the opposite direction.) 
  
We next introduce a special class of ``caterpillar'' %exchange  
patterns, and state conditions on their exchange polynomials  
that will imply Laurentness.

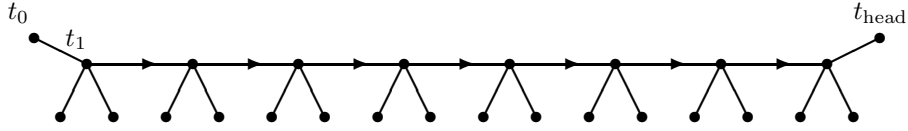
\begin{figure}[ht] 
\begin{center} 
\begin{picture}(280,60)(0,0) 
\thicklines 
 \put(0,30){\line(1,0){280}} 
 \multiput(20,30)(40,0){7}{\vector(1,0){7}} 
 \multiput(0,30)(40,0){8}{\line(1,-2){10}} 
 \multiput(0,30)(40,0){8}{\line(-1,-2){10}} 
 \put(0,30){\line(-2,1){20}} 
 \put(280,30){\line(2,1){20}} 
 \multiput(0,30)(40,0){8}{\circle*{4}} 
 \multiput(-10,10)(20,0){16}{\circle*{4}} 
 \put(-20,40){\circle*{4}} 
 \put(300,40){\circle*{4}} 
 \put(-30,47){$t_0$} 
 \put(-8,37){$t_1$} 
 \put(290,47){$t_{\rm head}$} 
\end{picture} 
\end{center} 
\caption{The ``caterpillar'' tree $\TT_{n,m}$, for $n=4$, $m=8$} 
\label{fig:caterpillar} 
\end{figure}

For $m \geq 1$, let $\TT_{n,m}$ be the tree of the form shown in 
Figure~\ref{fig:caterpillar}.  
 
The tree $\TT_{n,m}$ has $m$ vertices of degree $n$ in its 
``spine" and $m (n-2)+2$ vertices of degree~1. 
%; all vertices of degree $n$ lie on a path (``stem'') of length~$m$. 
We label every edge of the tree by an element of~$[n]$, so that the 
$n$ edges emanating from each vertex on the spine receive different 
labels. 
We let the root $t_0$ be a vertex in $\TT_{n,m}$ that does not belong 
to the spine but is connected to one of its ends. 
This gives rise to the orientation of the spine, with all the arrows 
pointing away from $t_0$ (see Figure~\ref{fig:caterpillar}). 
We assign a nonzero exchange polynomial $P\in\AAA[x_1,\dots,x_n]$  
to every edge $t\overunder{}{}t'$ of $\TT_{n,m}$,  
thus obtaining an exchange pattern. 
 
%%% To formulate the main result of this section, we will need some  
%%% preparation. 
%%% First, for 
For a rational function $F=F(x,y,\dots)$, 
we will denote by $F|_{x\leftarrow g(x,y,\dots)}$ 
the result of substituting $g(x,y,\dots)$ for $x$ into~$F$. 
To illustrate, if $F(x,y)=xy$, then $F|_{x\leftarrow 
 \frac{y}{x}}=\frac{y^2}{x}$.  
  
%%% We denote by $\AAA^{\times}$ the group of invertible elements  
%%% in~$\AAA$. 
%%% [Discuss divisibility and gcd]  
 
\begin{theorem} 
\label{th:laurent-gcd} 
{\rm (The Caterpillar Lemma)} 
Assume that a generalized exchange pattern on $\TT_{n,m}$  
satisfies the following conditions:  
\begin{eqnarray} 
\label{eq:GEP1} 
&&\text{For any edge $\bullet \overunder{k}{P} \bullet$,  
the polynomial $P$ does not depend on $x_k$, and is not }\\ 
\nonumber  
&&\text{divisible by any $x_i$, $i\in [n]$.}\\ 
\label{eq:GEP2} 
&&\text{For any two edges $\bullet \overunder{i}{P} \bullet 
\overunderarrow{j}{Q} \bullet$,  
the polynomials $P$ and $Q_0\!=\!Q|_{x_i=0}$}\\ 
\nonumber 
 && \text{are coprime elements of $\AAA[x_1,\dots,x_n]$.}\\ 
\label{eq:GEP3} 
&&\text{For any three edges 
$\bullet \overunder{i}{P} \bullet \overunderarrow{j}{Q} \bullet 
\overunder{i}{R} \bullet$  
labeled $i,j,i$, we have}\\ 
\nonumber  
%\label{eq:P=LR0/Q0} 
&& \hspace{1in} L\cdot Q_0^b \cdot P =R\bigl|_{x_j\leftarrow 
 \frac{Q_0}{x_j}},  
\\  
\nonumber 
&&\text{where $b$ is a nonnegative integer, $Q_0\!=\!Q|_{x_i=0}\,$,  
and $L$ is a Laurent}\\ 
\nonumber  
&& \text{monomial whose coefficient lies in $\AAA$ and 
is coprime with~$P$.}  
\end{eqnarray} 
Then each element $x_i(t)$, for $i\in [n]$, $t\in\TT_{n,m}\,$, 
is a Laurent polynomial 
in $x_1(t_0),\dots,x_n(t_0)$, with coefficients 
in~$\AAA$.  
\end{theorem} 
 
(Note the orientation of edges in (\ref{eq:GEP2})--(\ref{eq:GEP3}).)  
 
\proof 
Our argument is essentially the same as in \cite[Theorem~3.2]{fz-clust1}. 
For $t \in \TT_{n,m}$, let 
\[ 
\LL (t) = \AAA [x_1(t)^{\pm 1}, \dots, x_n(t)^{\pm 1}] 
\] 
denote the Laurent polynomial ring in the cluster $\textbf{x}(t)$ with coefficients in~$\AAA$. 
We view each $\LL (t)$ as a subring of the ambient field of rational  
functions~$\AAA(\textbf{x}(t_0))$. 
 
In this notation, our goal is to show that 
every cluster $\textbf{x}(t)$ is contained in $\LL (t_0)$. 
We abbreviate $\LL_0=\LL (t_0)$. 
Note that $\LL_0$ is a unique factorization domain, so 
any two elements $x, y \in \LL_0$ have a 
well-defined greatest common divisor $\gcd (x,y)$ 
which is an element of $\LL_0$ defined up to a multiple from the group  
$\LL_0^{\times}$ of invertible elements in~$\LL_0$;  
the group $\LL_0^{\times}$ consists of Laurent monomials  
in $x_1 (t_0), \dots, x_n (t_0)$ whose coefficient belongs to  
$\AAA^{\times}$, the group of invertible elements of $\AAA$.

To prove that all $\textbf{x}(t)$ are contained in $\LL_0$, 
we proceed by induction on $m$, the size of the spine. 
The claim is trivial for $m = 1$, so let us assume that 
$m \geq 2$, and furthermore assume that our statement is true for all 
``caterpillars" with smaller spine. 
It is thus enough to prove that $\textbf{x}(t_{\rm head})\subset  
\LL_0\,$, where $t_{\rm head}$ is one of the vertices most distant  
from $t_0$ (see Figure~\ref{fig:caterpillar}). 
 
We assume that the path from $t_0$ to $t_{\rm head}$ starts 
with the following two edges: 
$t_0 \overunder{i}{P} t_1 \overunderarrow{j}{Q} t_2$. 
Let $t_3 \in \TT_{n,m}$ be the vertex such that 
$t_2 \overunder{i}{R} t_3$.  
The following lemma plays a crucial role in our proof. 
 
\begin{lemma} 
\label{lem:3-step-general} 
%Let 
%$~~t_0 \overunder{i}{P} t_1 \overunderarrow{j}{Q} t_2 \overunder{i}{R} t_3~$ 
%be a fragment of a generalized exchange pattern. 
%Then $x_i (t_3) \in \LL_0=\LL(t_0)$, and therefore all variables in the 
The clusters ${\bf x}(t_1)$, ${\bf x}(t_2)$, and ${\bf x}(t_3)$ are contained 
in $\LL_0$. 
Furthermore, %we have 
$\gcd (x_i (t_3), x_i (t_1)) = \gcd (x_j (t_2), x_i (t_1)) = 1$.  
%(as elements of $\LL (t_0)$). 
\end{lemma}

\proof 
The only element in the clusters $\textbf{x}(t_1)$, $\textbf{x}(t_2)$, and 
$\textbf{x}(t_3)$ whose inclusion in $\LL_0$ is not 
immediate from (\ref{eq:exchange1})--(\ref{eq:exchange2}) 
is~$x_i (t_3)$. 
To simplify the notation, let us denote 
$x=x_i(t_0)$, $y=x_j(t_0)=x_j(t_1)$, $z=x_i(t_1)=x_i(t_2)$, 
$u=x_j(t_2)=x_j(t_3)$, and $v=x_i(t_3)$, so that these variables 
appear in the clusters at $t_0,\dots,t_3$, as shown below: 
\[ 
\begin{array}{c}\scriptstyle y,x\\ \bullet\\ \scriptstyle t_0\end{array} 
 \hspace{-.1in}\wideoverunder{i}{P}\hspace{-.1in} 
\begin{array}{c}\scriptstyle z,y\\ \bullet\\ \scriptstyle t_1\end{array} 
 \hspace{-.1in}\wideoverunderarrow{j}{Q}\hspace{-.1in} 
\begin{array}{c}\scriptstyle u,z\\ \bullet\\ \scriptstyle t_2\end{array} 
 \hspace{-.1in}\wideoverunder{i}{R}\hspace{-.1in} 
\begin{array}{c}\scriptstyle v,u\\ \bullet\\ \scriptstyle t_3\end{array} 
\,. 
\] 
Note that the variables $x_k$, for $k\notin\{i,j\}$, do not change as 
we move among the four clusters under consideration. 
The lemma is then restated as saying that 
\begin{eqnarray} 
\label{eq:Laurent-3} 
&&\text{$v\in\LL_0$;}\\ 
\label{eq:gcd(u,z)} 
&&\text{$\gcd(z,u)=1$ 
%(as elements of $\LL_0$)  
;}\\ 
\label{eq:gcd(v,z)} 
&&\text{$\gcd(z,v)=1$ 
%(as elements of $\LL_0$)  
.} 
\end{eqnarray} 
Another notational convention will be based on the fact that each of 
the polynomials $P,Q,R$ has a distinguished variable on which it 
depends, namely $x_j$ for $P$ and $R$, and $x_i$ for~$Q$. 
(In view of (\ref{eq:GEP1}), $P$ and $R$ do not depend on $x_i$, while 
$Q$ does not depend on~$x_j$.) 
With this in mind, we will routinely write $P$, $Q$, and $R$ as 
polynomials in one (distinguished) variable. 
For example, we rewrite the formula in (\ref{eq:GEP3}) as  
\begin{equation} 
\label{eq:R(Q(0)/y)=...P} 
R\left(\frac{Q(0)}{y}\right)=L(y) Q(0)^b P(y), 
\end{equation} 
where we denote 
$L(y)=L|_{x_j\leftarrow y}$.  
In the same spirit, the notation $Q'$, $R'$, etc., will refer to the 
partial derivatives with respect to the distinguished variable. 
 
We will prove the statements (\ref{eq:Laurent-3}), 
(\ref{eq:gcd(u,z)}), and (\ref{eq:gcd(v,z)}) one by one, in this 
order. 
We have: 
\begin{eqnarray*} 
&&z=\frac{P(y)}{x}\,;\\ 
&&u=\frac{Q(z)}{y}=\frac{Q\left(\frac{P(y)}{x}\right)}{y}\,;\\ 
&&v=\frac{R(u)}{z} 
=\frac{R\left(\frac{Q(z)}{y}\right)}{z} 
=\frac{R\left(\frac{Q(z)}{y}\right)-R\left(\frac{Q(0)}{y}\right)}{z} 
+\frac{R\left(\frac{Q(0)}{y}\right)}{z}\,. 
\end{eqnarray*} 
Since 
\[ 
\frac{R\left(\frac{Q(z)}{y}\right)-R\left(\frac{Q(0)}{y}\right)}{z} 
\in\LL_0 
\] 
and 
\[ 
\frac{R\left(\frac{Q(0)}{y}\right)}{z} 
= \frac{L(y) Q(0)^b P(y)}{z} 
= L(y) Q(0)^b x  
\in\LL_0 \,, 
\] 
(\ref{eq:Laurent-3}) follows. 
 
We next prove (\ref{eq:gcd(u,z)}). 
We have 
\[ 
u=\frac{Q(z)}{y}\equiv \frac{Q(0)}{y}\bmod z\,. 
\] 
Since $x$ and $y$ are invertible in $\LL_0$, we conclude that 
$\gcd(z,u)=\gcd(P(y),Q(0))=1$ (using (\ref{eq:GEP2})). 
 
It remains to prove (\ref{eq:gcd(v,z)}). 
Let 
\[ 
f(z)=R\left(\textstyle\frac{Q(z)}{y}\right). 
\] 
Then 
\[ 
v=\frac{f(z)-f(0)}{z}+L(y) Q(0)^b x\,. 
\] 
Working $\bmod z$, we obtain:  
\[ 
\frac{f(z)-f(0)}{z} 
\equiv f'(0) 
=R'\left(\textstyle\frac{Q(0)}{y}\right) 
\cdot\textstyle\frac{Q'(0)}{y}\,. 
\] 
%%% 
%%% From (\ref{eq:R(u)=A(u)P}) we get 
%%% \[ 
%%% R'(u) 
%%% =A'(u) \, P\!\left(\textstyle\frac{Q(0)}{u}\right) 
%%% -A(u)\, \frac{Q(0)}{u^2}\, P'\!\left(\textstyle\frac{Q(0)}{u}\right) , 
%%% \] 
%%% implying 
%%% \begin{eqnarray*} 
%%% R'\!\left(\textstyle\frac{Q(0)}{y}\right) 
%%% &=&A'\!\left(\textstyle\frac{Q(0)}{y}\right) \, P(y) 
%%% -A\!\left(\textstyle\frac{Q(0)}{y}\right)\, \frac{y^2}{Q(0)} 
%%% \, P'(y) \\ 
%%% &\equiv& -A\!\left(\textstyle\frac{Q(0)}{y}\right)\, \frac{y^2}{Q(0)} 
%%% \, P'(y)~ \bmod z\,. 
%%% \end{eqnarray*} 
%%% %(since $P(y)\equiv 0\bmod z$). 
%%% We conclude that, $\bmod z$, 
%%% \[ 
%%% v\equiv -A\!\left(\textstyle\frac{Q(0)}{y}\right)\, \frac{y Q'(0)}{Q(0)} 
%%% \, P'(y) 
%%% +A\left(\textstyle\frac{Q(0)}{y}\right)x\,. 
%%% \] 
Hence 
\[ 
v\equiv R'\left(\textstyle\frac{Q(0)}{y}\right) 
\cdot\textstyle\frac{Q'(0)}{y} 
+L(y) Q(0)^b x 
\bmod z\,. 
\] 
Note that the right-hand side is a polynomial of degree~1 in $x$ 
whose coefficients are Laurent polynomials in the rest of the 
variables of the cluster $\textbf{x}(t_0)$. 
Thus (\ref{eq:gcd(v,z)}) follows from  
$\gcd\left(L(y) Q(0)^b,P(y)\right)=1$, 
which is a consequence of (\ref{eq:GEP2})--(\ref{eq:GEP3}).  
\endproof

We can now complete the proof of Theorem~\ref{th:laurent-gcd}. 
We need to show that %for any vertices $t_0,t\in T$, 
any variable $X= x_k(t_{\rm head})$ %, for $t\in T$ and $k\in[n]$, 
belongs to $\LL_0$. 
Since both $t_1$ and $t_3$ are closer to $t_{\rm head}$ than 
$t_0$, we can use the inductive assumption to conclude 
that $X$ belongs to both $\LL(t_1)$ and $\LL(t_3)$. 
Since $X \in \LL(t_1)$, it follows from (\ref{eq:exchange1}) that 
$X$ can be written as $X = f/ x_i(t_1)^a$ for some $f \in \LL_0$ 
and $a\in\ZZ_{\geq 0}\,$. 
On the other hand, since $X \in \LL(t_3)$, it follows from 
(\ref{eq:exchange1}) and from the inclusion $x_i (t_3) \in \LL_0$ 
provided by Lemma~\ref{lem:3-step-general} 
that $X$ has the form $X = g/ x_j(t_2)^b x_i(t_3)^c$ for some $g \in \LL_0$ 
and some $b, c\in\ZZ_{\geq 0}\,$. 
The inclusion $X \in \LL_0$ now follows from the fact 
that, by the last statement in Lemma~\ref{lem:3-step-general}, 
the denominators in the two obtained expressions for $X$ are 
coprime in $\LL_0$. 
\endproof

\section{One-dimensional recurrences} 
\label{sec:gale-robinson} 
 
In this section, we apply Theorem~\ref{th:laurent-gcd} 
to study the Laurent phenomenon for sequences 
$y_0,y_1,\dots$ given by recursions of the form 
\begin{equation} 
\label{eq:one-dim-recursion} 
y_{m+n}y_m = F(y_{m+1},\dots,y_{m+n-1}), 
\end{equation} 
where $F\in\AAA[x_1,\dots,x_{n-1}]$. 
 
For an integer $m$, let $\rem{m}$ denote the unique element of 
$[n]=\{1,\dots,n\}$ satisfying $m\equiv\rem{m} \bmod n$. 
We define the polynomials $F_1,\dots,F_n\in\AAA[x_1,\dots,x_n]$ by 
\begin{equation} 
\label{eq:Fm-cyclic} 
F_m=F(x_{\rem{m+1}},x_{\rem{m+2}},\dots,x_{\rem{m-1}}); 
\end{equation} 
thus $F_m$ does not depend on the variable~$x_m$. 
We introduce the infinite ``cyclic exchange pattern'' 
\begin{equation} 
\label{eq:cyclic1} 
%\cdots \overunder{}{} 
%\bullet\wideoverunder{2}{a_{-1}x_1+b_{-1}} 
%t_0 \wideoverunder{1}{} 
t_0 \wideoverunder{\rem{0}}{F_{\rem{0}}} 
t_1 \wideoverunder{\rem{1}}{F_{\rem{1}}} 
t_2 \wideoverunder{\rem{2}}{F_{\rem{2}}} 
t_3 \wideoverunder{\rem{3}}{F_{\rem{3}}} 
t_4 \overunder{}{} 
\cdots \,, 
\end{equation} 
and let the cluster at each point $t_m$ consist of the variables 
$y_m,\dots,y_{m+n-1}$, labeled within the cluster according to the rule 
$y_s=x_{\rem{s}}(t_m)$. 
Then equations (\ref{eq:one-dim-recursion}) become 
 the exchange relations associated with this pattern. 
 
To illustrate, let $n=4$. 
Then the clusters will look like this: 
\[ 
%\cdots \overunder{1}{}\hspace{-.2in} 
\begin{array}{c}\scriptstyle y_1,y_2,y_3,y_0\\ \bullet\\ t_0\end{array} 
 \hspace{-.2in}\wideoverunder{4}{}\hspace{-.2in} 
\begin{array}{c}\scriptstyle y_1,y_2,y_3,y_4\\ \bullet\\ t_1\end{array} 
 \hspace{-.2in}\wideoverunder{1}{}\hspace{-.2in} 
\begin{array}{c}\scriptstyle y_5,y_2,y_3,y_4\\ \bullet\\ t_2\end{array} 
 \hspace{-.2in}\wideoverunder{2}{}\hspace{-.2in} 
\begin{array}{c}\scriptstyle y_5,y_6,y_3,y_4\\ \bullet\\ t_3\end{array} 
 \hspace{-.2in}\wideoverunder{3}{}\hspace{-.2in} 
\begin{array}{c}\scriptstyle y_5,y_6,y_7,y_4\\ \bullet\\ t_4\end{array} 
 \hspace{-.2in}\wideoverunder{4}{} 
\cdots \,. 
\] 
In order to include this situation into the setup of 
Section~\ref{sec:caterpillar} (cf.\ Figure~\ref{fig:caterpillar}), 
we create an infinite ``caterpillar tree'' whose ``spine'' is formed  
by the vertices $t_m$, $m>0$.  
We thus attach the missing $n-2$ ``legs'' with labels in  
$[n]-\{\rem{m-1},\rem{m}\}$, to each vertex~$t_m$. 
  
Our next goal is to state conditions on the polynomial~$F$  
which make it possible to assign exchange polynomials satisfying  
(\ref{eq:GEP1})--(\ref{eq:GEP3}) to the newly constructed legs. 
The first requirement (cf.\ (\ref{eq:GEP1})) is:  
\begin{equation}  
\label{eq:GEP1a}  
\text{The polynomial $F$ is not divisible by any $x_i$, $i\in [n-1]$.}  
\end{equation}  
For $m\in[n-1]$, we set  
\begin{equation}  
\label{eq:Qm}  
Q_m=F_{m}|_{x_n\leftarrow 0}  
=F(x_{m+1},\dots,x_{n-1},0,x_1,\dots,x_{m-1}). 
\end{equation}  
Our second requirement is  
\begin{equation}  
\label{eq:GEP2a}  
\text{Each $Q_m$ %, for $m\in[n-1]$, 
is an irreducible element of 
$\AAA[x_1^{\pm1},\dots,x_{n-1}^{\pm1}]$.}  
\end{equation}  
  
\begin{figure}
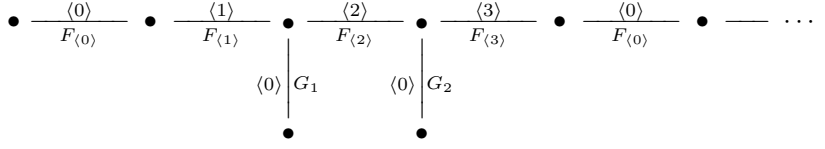
 
\[ 
\bullet 
\wideoverunder{\rem{0}}{F_{\rem{0}}} 
\bullet 
\wideoverunder{\rem{1}}{F_{\rem{1}}}\hspace{-.2in} 
\begin{array}{c}\\[.43in] \bullet\\ {}_{\rem{0}} \Biggl| {}_{G_1}\\ 
\bullet\end{array} 
 \hspace{-.2in}\wideoverunder{\rem{2}}{F_{\rem{2}}} 
\hspace{-.2in} 
\begin{array}{c}\\[.43in] \bullet\\ {}_{\rem{0}} \Biggl| {}_{G_2}\\ 
\bullet\end{array} 
 \hspace{-.2in} 
\wideoverunder{\rem{3}}{F_{\rem{3}}} 
\bullet 
\wideoverunder{\rem{0}}{F_{\rem{0}}} 
\bullet 
\overunder{}{} 
\cdots \,. 
\] 
\caption{Constructing a caterpillar; $n=4$.} 
\label{fig:cater4} 
\end{figure}

To state our most substantial requirement, we recursively define a sequence of  
polynomials $G_{n-1},\dots,G_1,G_0\in\AAA[x_1,\dots,x_{n-1}]$;  
more precisely, each $G_m$ will be defined up to a multiple in~$\AAA^\times$. 
(Later, $G_1,\dots,G_{n-2}$ will become the exchange polynomials assigned  
to the ``legs'' of the caterpillar labeled by~$n=\rem{0}$; see 
Figure~\ref{fig:cater4}.) 
  
We set $G_{n-1}=F$, and obtain each $G_{m-1}$ from $G_m$, as follows.  
Let  
\begin{equation}  
\label{eq:G-step1}  
\stackrel{\sim}{G}_{m-1} = G_m\bigl|_{x_m\leftarrow \frac{Q_m}{x_m}}\ .  
\end{equation}  
Let $L$ be a Laurent monomial in $x_1,\dots,x_{n-1}$, with coefficient in  
$\AAA$, such that 
\begin{equation}  
\label{eq:G-step2}  
\stackrel{\approx}{G}_{m-1}  
=\frac{\stackrel{\sim}{G}_{m-1}}{L}  
\end{equation}  
is a polynomial in $\AAA[x_1,\dots,x_{n-1}]$  
not divisible by any $x_i$ or by any non-invertible scalar  
in~$\AAA$. 
Such an $L$ is unique up to a multiple in~$A^\times$. 
Finally, we set  
\begin{equation}  
\label{eq:G-step3}  
G_{m-1}=\frac{\stackrel{\approx}{G}_{m-1}}{Q_m^b}\,,  
\end{equation}  
where $Q_m^b$ is the maximal power of $Q_m$ that divides 
$\stackrel{\approx}{G}_{m-1}$. 
With all this notation, our final requirement is: 
\begin{equation}  
\label{eq:GEP3a}  
G_0=F. %\,,\text{~where~} \alpha\in\AAA^\times. 
\end{equation}  
  
%%% Let us formally describe this recursive algorithm. 
%%%  
%%% \vbox{\begin{algorithm} 
%%% \label{alg:cyclic} 
%%% {\rm \emph{Input:} Polynomial $F\in\AAA[x_1,\dots,x_{n-1}]$. 
%%% \emph{Output:} Polynomial~$G$. 
%%%  
%%% \begin{itemize} 
%%% \item[\textbf{1.}] Set $m:=n-1$; $G:=F$. 
%%%  
%%% \item[\textbf{2.}] Set 
%%% %\begin{equation} 
%%% %\label{eq:Qm} 
%%% $Q:=Q_m:=F(x_{m+1},\dots,x_{n-1},0,x_1,\dots,x_{m-1})$.\\ 
%%% (\emph{Comment:} $Q_m=F_{\rem{m}}|_{x_n\leftarrow 0}$.) 
%%%  
%%% \item[\textbf{3.}] Set $G:=G\bigl|_{x_m\leftarrow \frac{Q}{x_m}}$. 
%%%  
%%% \item[\textbf{4.}] Divide $G$ by the maximal power of $Q$ that divides it. 
%%%  
%%% \item[\textbf{5.}] Multiply $G$ by a Laurent monomial so that it became a 
%%% polynomial not divisible by any~$x_i$. 
%%% (\emph{Comment:} the result is defined up to a factor in~$\AAA^{\times}$.) 
%%%  
%%% \item[\textbf{6.}] Set $m:=m-1$. (\emph{Comment:} the current  
%%% value of $G$ is~$G_m$.)  
%%%  
%%% \item[\textbf{7.}] If $m>0$, go to step~\textbf{2}. 
%%% Otherwise stop. 
%%% \end{itemize} 
%%% } 
%%% \end{algorithm} 
%%% } 

\begin{theorem} 
\label{th:laurent-cyclic} 
Let $F$ be a polynomial in the variables $x_1,\dots,x_{n-1}$ with 
coefficients in a unique factorization domain~$\AAA$  
satisfying conditions {\rm (\ref{eq:GEP1a}), (\ref{eq:GEP2a}),  
and (\ref{eq:GEP3a})}. 
Then every term of the sequence $(y_i)$ defined by the 
recurrence 
\begin{equation*} 
%\label{eq:laurent-cyclic} 
y_{m+n} = \frac{F(y_{m+1},\dots,y_{m+n-1})}{y_m} 
\end{equation*} 
%{\rm(\ref{eq:one-dim-recursion})} 
is a Laurent polynomial in the initial $n$ terms, 
with coefficients in~$\AAA$. 
\end{theorem} 
  
\proof  
%%% The proof reduces to a careful verification of the conditions in  
%%% Theorem~\ref{th:laurent-gcd} for a caterpillar pattern defined by a  
%%% right-to-left recursive procedure sketched above. 
%%% The exchange polynomials assigned to the legs of the pattern will be  
%%% cyclic shifts of the polynomials $G_m$, which were devised with the  
%%% express purpose of satisfying {\rm  
%%% (\ref{eq:GEP1})--(\ref{eq:GEP3})}. 
%%% Now, let us make this precise. 
To prove the Laurentness of some $y_N$,  
we will apply Theorem~\ref{th:laurent-gcd} 
to the caterpillar tree constructed as follows. 
We set $t_{\rm head}=t_{N-n+1}$;  
this corresponds to the first cluster containing~$y_N$. 
As a path from $t_0$ to $t_{\rm head}$, 
we take a finite segment of (\ref{eq:cyclic1}): 
\begin{equation} 
\label{eq:t0-tN-n}  
t_0 \wideoverunder{\rem{0}}{F_{\rem{0}}} 
t_1 \wideoverunder{\rem{1}}{F_{\rem{1}}} 
t_2 \wideoverunder{\rem{2}}{F_{\rem{2}}} 
\cdots  
\wideoverunder{\rem{N-1}}{F_{\rem{N-1}}} 
t_{N-n}  
\wideoverunder{\rem{N}}{F_{\rem{N}}} 
t_{N-n+1}\,. 
\end{equation} 
We then define the exchange polynomial $G_{j,k-1}$ 
associated with the leg labeled $j$ attached  
to a vertex $t_k$ on the spine (see Figure~\ref{fig:leg})  
by 
\[  
G_{j,k-1}=G_{\rem{k-j-1}}(x_{\rem{j+1}},\dots,x_n,x_1,\dots,x_{\rem{j-1}}),  
\]  
where in the right-hand side, we use the polynomials  
$G_1,\dots,G_{n-2}$ constructed in  
(\ref{eq:G-step1})--(\ref{eq:G-step3}) above. 
  
\begin{figure}[ht]  
\vspace{-.5in}  
\begin{center}  
$  
\bullet  
\wideoverunder{\rem{k-1}}{F_{\rem{k-1}}} 
\hspace{-.35in} 
\begin{array}{c}\\[.43in] t_k\\ \hspace{.28in} {\scriptstyle j} \Biggl| {\scriptstyle G_{j,k-1}}\\ 
\bullet\end{array} 
 \hspace{-.4in} 
\wideoverunder{\rem{k}}{F_{\rem{k}}} 
\bullet 
$  
\end{center}  
\caption{}  
\label{fig:leg}  
\end{figure}  
  
It remains to verify that this exchange pattern satisfies  
(\ref{eq:GEP1}), (\ref{eq:GEP2}), and (\ref{eq:GEP3}). 
Condition (\ref{eq:GEP1}) for the edges appearing in  
(\ref{eq:t0-tN-n}) is immediate from (\ref{eq:GEP1a}),  
while for the rest of the edges, it follows from the definition of  
$\stackrel{\approx}{G}_{m-1}\,$ in (\ref{eq:G-step2}). 
  
Turning to (\ref{eq:GEP2}), we first note that we may assume  
$i=\rem{0}=n$ (otherwise apply a cyclic shift of indices). 
Under this assumption, we can identify the polynomials $P$ and $Q_0$ in  
(\ref{eq:GEP2}) with the polynomials $G_{m-1}$ and $Q_m$ in  
(\ref{eq:G-step3}), for some value of~$m$. 
(The special case of $P$ attached to one of the edges in  
(\ref{eq:t0-tN-n}) corresponds to $m=1$, and its validity requires  
(\ref{eq:GEP3a}).) 
Then the condition $\gcd(G_{m-1},Q_m)=1$ follows from (\ref{eq:GEP2a})  
and the choice of the exponent $b$ in (\ref{eq:G-step3}). 
  
Finally, (\ref{eq:GEP3}) is ensured by the construction 
(\ref{eq:G-step1})--(\ref{eq:G-step3}), which was designed expressly for  
this purpose. As before, the special case 
of $P$ attached to one of the edges in (\ref{eq:t0-tN-n}) holds due to 
(\ref{eq:GEP3a}). 
\endproof  
 
%%% In this diagram, the polynomial $G_2$ is constructed from 
%%% $F_{\rem{0}}$ and $F_{\rem{3}}$ using~(c); 
%%% then $G_1$ is constructed from 
%%% $G_2$ and $F_{\rem{2}}$, again using~(c); 
%%% finally, we apply the same construction to $G_1$ and $F_{\rem{1}}$, 
%%% and if we are lucky enough to recover $F_{\rem{0}}$, 
%%% then Theorem~\ref{th:laurent-gcd} applies, and the Laurent phenomenon 
%%% follows. 
%%% Note that we will not have to repeat the argument for the legs 
%%% carrying labels other than $\rem{0}$, since those cases can be 
%%% obtained from the $\rem{0}$ case by a shift of index. 
%%% (Remember that the polynomials $F_{\rem{m}}$ differ from each other by 
%%% cyclic permutation of the variables.) 
  
In the rest of this section, we give a few applications of  
Theorem~\ref{th:laurent-cyclic}. 
In all of them, conditions (\ref{eq:GEP1a}) and  
(\ref{eq:GEP2a}) are immediate, so we concentrate on the verification  
of (\ref{eq:GEP3a}). 
 
\begin{example} 
\label{example:somos3}  
{\rm 
Let $a$ and $b$ be positive integers, 
and let the sequence $y_0,y_1,\dots$ satisfy the recurrence 
\begin{equation*} 
y_{k}=\frac{y_{k-2}^{a}y_{k-1}^{b}+1}{y_{k-3}}\,. 
\end{equation*} 
%(This is a generalization of the ordinary Somos-3 sequence, 
%which arises in the special case $a=b=1$.) 
We claim that every term of the sequence is a Laurent polynomial over 
$\ZZ$ in $y_0$, $y_1$, and~$y_2$. 
%(This example can be further generalized by adding coefficients, 
%and by making the exponents dependent on $k\bmod 3$.) 
To prove this, we set $n=3$ and construct the polynomials  
$G_2$, $G_1$, and $G_0$ using (\ref{eq:G-step1})--(\ref{eq:G-step3}). 
Initializing $G_2=F(x_1,x_2)=x_1^a x_2^b+1$, we obtain: 
\[ 
\begin{array}{lll} 
Q_2=F(0,x_1)=1,\quad 
& \stackrel{\sim}{G}_{1}=F\bigl|_{x_2\leftarrow\frac{Q_2}{x_2}} 
=x_1^a x_2^{-b}+1, \quad 
& G_1=\stackrel{\approx}{G}_{1}=x_1^a+x_2^b,\\[.1in] 
Q_1=F(x_2,0)=1, 
& \stackrel{\sim}{G}_{0}=G_1\bigl|_{x_1\leftarrow\frac{Q_1}{x_1}} 
=x_1^{-a}+x_2^b, \quad 
& G_0=\stackrel{\approx}{G}_{0}=1+x_1^a x_2^b=F, 
\end{array} 
\] 
as desired. 
} 
\end{example} 
 
\begin{example} 
\label{example:somos4}  
(Generalized Somos-4 sequence) 
{\rm 
Let $a$, $b$, and $c$ be positive integers, 
and let the sequence $y_0,y_1,\dots$ satisfy the recurrence 
\begin{equation*} 
y_{k}=\frac{y_{k-3}^{a}y_{k-1}^{c}+y_{k-2}^b}{y_{k-4}}\,. 
\end{equation*} 
(The Somos-4 sequence~\cite{gale-intelligencer}, introduced by Michael Somos,  
is the special case $a=c=1$, $b=2$.) 
Again, each $y_i$ is a Laurent polynomial in the initial terms  
$y_0$, $y_1$, $y_2$, and~$y_3$. 
To prove this, we set $n=4$ and compute 
$G_3,\dots,G_0$ using (\ref{eq:G-step1})--(\ref{eq:G-step3}) and 
beginning with  
$G_3=F=x_1^a x_3^c+x_2^b$: 
\[ 
\begin{array}{lll} 
\!\!Q_3\!=\!F(0,x_1,x_2)\!=\!x_1^b, 
 & G_3\bigl|_{x_3\leftarrow\frac{Q_3}{x_3}} 
 \!=\!x_1^{a+bc}x_3^{-c}+x_2^b, 
   & G_2\!=\!x_1^{a+bc} +x_2^b x_3^c,\\[.1in] 
\!\!Q_2\!=\!F(x_3,0,x_1)\!=\!x_1^c x_3^a,~~ 
 & G_2\bigl|_{x_2\leftarrow\frac{Q_2}{x_2}} 
 \!=\!x_1^{a+bc} \!+\!x_1^{bc}x_2^{-b} x_3^{ab+c}\!, ~~ 
   & G_1\!=\!x_1^a x_2^b + x_3^{ab+c},\\[.1in] 
\!\!Q_1\!=\!F(x_2,x_3,0)\!=\!x_3^b, 
 & G_1\bigl|_{x_1\leftarrow\frac{Q_1}{x_1}} 
 \!=\!x_1^{-a} x_2^b x_3^{ab} + x_3^{ab+c}, 
   & G_0\!=\!x_2^b \!+\! x_1^a x_3^c \!=\!F, 
\end{array} 
\] 
and the claim follows. 
} 
\end{example} 

\begin{remark}
\label{rem:GR-stuff}
{\rm The Laurent phenomena in Theorems~\ref{th:gale-robinson} 
and \ref{th:two-term-gale-robinson} 
can also be proved by applying Theorem~\ref{th:laurent-cyclic}: 
in the former (resp., latter) case, the polynomial $F$ is given by
$F=\alpha x_{p}x_{n-p}+\beta x_{q}x_{n-q} +\gamma x_{r}x_{n-r}$
(resp., 
$F=\alpha x_{p}x_{n-p}+\beta x_{q}x_{n-q}$).
The proofs are straightforward but rather long.
Shorter proofs, based on J.~Propp's idea of viewing one-dimensional
recurrences as ``projections'' of multi-dimensional ones, 
are given in Section~\ref{sec:lattice-based} below.
}
\end{remark}

\section{%Number walls, frieze patterns, and other  
Two- and three-dimensional recurrences} 
\label{sec:lattice-based} 
 
In this section, we use the strategy of  
Section~\ref{sec:gale-robinson} 
% Theorem~\ref{th:laurent-gcd} 
to establish the Laurent phenomenon for several recurrences  
involving two- and three-dimensional arrays. 
Our first example generalizes a construction (and the corresponding Laurentness  
conjecture) suggested by Noam Elkies and communicated by James Propp. 
Even though the Laurent phenomenon in this example can be deduced 
from Theorem~\ref{th:oct}, 
we choose to give a self-contained treatment, for the sake of exposition. 
%We will discuss this example in considerable detail, for it will serve as a 
%prototype for the subsequent applications. 
%Those will not be treated nearly as meticulously, 
%as we will trust the reader with filling out the details.  
 
\begin{example} %(Elkies' recurrence) 
\label{example:elkies} 
(The knight recurrence) 
{\rm 
Consider a two-dimensional array  
$(y_{ij})_{i,j\geq 0}$ 
whose entries satisfy the recurrence 
\begin{equation} 
\label{eq:elkies} 
y_{i,j}y_{i-2,j-1}=\alpha y_{i,j-1}y_{i-2,j} +\beta y_{i-1,j} y_{i-1,j-1}\,.  
\end{equation} 
We will prove that every $y_{ij}$ is a 
Laurent polynomial in the initial entries  
\[ 
Y_{\rm init} = \{y_{ab}\,:\,a<2\text{~or~} b<1\},  
\]  
with coefficients in the ring $\AAA=\ZZ[\alpha,\beta]$. 
} 
\end{example} 
  
We will refer to $Y_{\rm init}$ as the \emph{initial cluster,}  
even though it is an infinite set.  
Notice, however, that each individual $y_{ij}$ only depends on  
finitely many variables $\{y_{ab}\in Y_{\rm init}\,:\,a\leq i, \ b\leq j\}$. 
  
Similarly to Section~\ref{sec:gale-robinson}, we will use the exchange 
relations (\ref{eq:elkies}) to create a sequence of 
clusters satisfying the Caterpillar Lemma 
(Theorem~\ref{th:laurent-gcd}). 
  
\begin{figure}[ht] 
\setlength{\unitlength}{4pt} 
 
\centering 
\begin{picture}(60,35)(0,7) 
\linethickness{0.05pt}  
\multiput(0,15)(0,5){5}{\line(1,0){60}} 
\multiput(0,10)(5,0){13}{\line(0,1){30}} 
\multiput(0,10)(0,30){2}{\line(1,0){60}} 
\thicklines  
\put(0,10){\line(2,1){60}}  
\put(0,15){\line(2,1){50}}  
\put(0,20){\line(2,1){40}}  
\put(0,25){\line(2,1){30}}  
\put(0,30){\line(2,1){20}}  
\put(0,35){\line(2,1){10}}  
  
\put(5,10){\line(2,1){55}}  
\put(5,15){\line(2,1){50}}  
\put(5,20){\line(2,1){40}}  
\put(5,25){\line(2,1){30}}  
\put(5,30){\line(2,1){20}}  
\put(5,35){\line(2,1){10}}  
  
\put(10,10){\line(2,1){50}}  
\put(15,10){\line(2,1){45}}  
\put(20,10){\line(2,1){40}}  
\put(25,10){\line(2,1){35}}  
\put(30,10){\line(2,1){30}}  
\put(35,10){\line(2,1){25}}  
\put(40,10){\line(2,1){20}}  
\put(45,10){\line(2,1){15}}  
\put(50,10){\line(2,1){10}}  
\put(55,10){\line(2,1){5}}  
  
%\thinlines 
\multiput(10,15)(5,0){11}{\circle{1}} 
\multiput(10,20)(5,0){11}{\circle{1}} 
\multiput(10,25)(5,0){11}{\circle{1}} 
\multiput(10,30)(5,0){11}{\circle{1}} 
\multiput(10,35)(5,0){11}{\circle{1}} 
\multiput(10,40)(5,0){11}{\circle{1}} 
  
\multiput(0,10)(5,0){13}{\circle*{1}} 
\multiput(0,15)(0,5){6}{\circle*{1}} 
\multiput(5,15)(0,5){6}{\circle*{1}} 
\put(0,7){\makebox(0,0){$0$}} 
%\put(60,7){\makebox(0,0){$N$}} 
\put(64,30){\makebox(0,0){$\mathcal{H}$}} 
\put(-3,10){\makebox(0,0){$0$}} 
%\put(-3,40){\makebox(0,0){$M$}} 

\end{picture} 
 
\caption{The initial cluster and the equivalence classes $\rem{h}$} 
\label{fig:elkies} 
\end{figure}
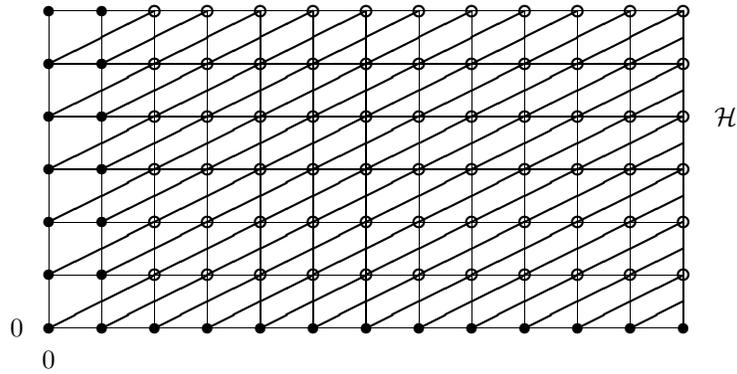

This is done in the following way.  
Let us denote by 
$\mathcal{H}=\ZZ_{\geq 0}^2$ the underlying set of indices; 
for $h=(i,j)\in \mathcal{H}$, we will write $y_h=y_{ij}\,$.  
The variables of the initial cluster have labels in the set  
\[ 
H_{\rm init} = \{(i,j)\in\mathcal{H}\,:\,i<2 \text{~or~} j<1\}.  
\] 
In Figure~\ref{fig:elkies}, the elements of $H_{\rm init}$ are marked 
by~$\bullet$'s.  
 
We introduce the product partial order on $\mathcal{H}$: 
\begin{equation}  
\label{eq:product-order} 
(i_1,j_1)\leq (i_2,j_2) \stackrel{\rm def}{\Leftrightarrow} 
(i_1\leq i_2)\text{~and~} (j_1\leq j_2).  
\end{equation}  
For an element $h=(i,j)\in \mathcal{H} - H_{\rm init}\,$,  
let us denote $h^-=(i-2,j-1)$; 
in this notation, the exchange relation (\ref{eq:elkies})  
expresses the product $y_h\cdot y_{h^-}$ as a polynomial in the 
variables $y_{h'}\,$, for $h^-< h'< h$. 
 
We write $h^-\sim h$, and extend this to an equivalence relation  
$\sim$ on $\mathcal{H}$. The equivalence class of $h$ is denoted  
by~$\rem{h}$. 
These classes are shown as slanted lines in Figure~\ref{fig:elkies}. 
All our exchange polynomials will belong to the ring 
$\AAA[x_a : a \in\mathcal{H}/\!\!\sim\,]$.  
  
\begin{figure}[ht] 
\setlength{\unitlength}{4pt} 
 
\centering 
  \begin{picture}(60,35)(0,7) 
%\thicklines 
\multiput(0,15)(0,5){5}{\line(1,0){60}} 
\multiput(0,10)(5,0){13}{\line(0,1){30}} 
\multiput(0,10)(0,30){2}{\line(1,0){60}} 
%\thinlines 
\multiput(15,10)(5,0){10}{\circle*{1}} 
\multiput(0,30)(0,5){3}{\circle*{1}} 
\multiput(5,25)(0,5){4}{\circle*{1}} 
\multiput(10,20)(0,5){3}{\circle*{1}} 
\multiput(15,15)(0,5){3}{\circle*{1}} 
\multiput(20,15)(0,5){2}{\circle*{1}} 
\put(0,7){\makebox(0,0){$0$}} 
%\put(60,7){\makebox(0,0){$M$}} 
%\put(64,30){\makebox(0,0){$\mathcal{H}$}} 
\put(-3,10){\makebox(0,0){$0$}} 
%\put(-3,40){\makebox(0,0){$N$}} 
\end{picture} 
 
\caption{Indexing set $H_9$} 
\label{fig:elkies-H9} 
\end{figure}
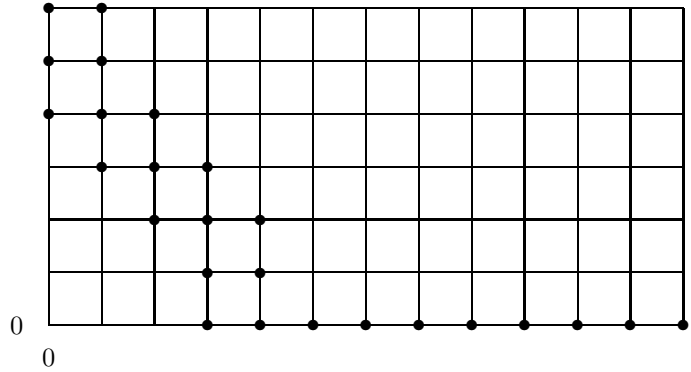

Note that $H_{\rm init}$ has exactly one representative from each  
equivalence class. 
We will now construct a sequence of subsets  
$H_0=H_{\rm init}, H_1,H_2,\dots$, each having this property, 
using the following recursive rule. 
Let us fix a particular linear extension of the partial order  
(\ref{eq:product-order}), say, 
\[ 
(i_1,j_1)\preceq (i_2,j_2) \stackrel{\rm def}{\Leftrightarrow} 
(i_1+j_1<i_2+j_2)\text{~or~} (i_1+j_1= i_2+j_2 \text{~and~} i_1\leq i_2). 
\] 
Restricting this linear ordering to the complement $\mathcal{H}-H_{\rm  
init}$ of the initial cluster, we obtain a numbering of the elements  
of this complement by positive integers:  
\[  
\begin{array}{l}  
h_0=(2,1),~h_1=(2,2),~h_2=(3,1),~h_3=(2,3),~h_4=(3,2),\\  
h_5=(4,1),~h_6=(2,4),~h_7=(3,3),~h_8=(4,2),  
\end{array}  
\]  
and so on. 
Having constructed $H_{m}$, we let $H_{m+1}=H_{m}\cup\{h_m\}-\{h_m^-\}$.  
To illustrate, the set $H_9$ is shown in Figure~\ref{fig:elkies-H9}.  
 
We next create the infinite exchange pattern  
\begin{equation} 
\label{eq:knight-spine} 
t_0 \wideoverunder{\rem{h_0}}{P_{\rem{h_0}}} 
t_1 \wideoverunder{\rem{h_1}}{P_{\rem{h_1}}} 
t_2 \wideoverunder{\rem{h_2}}{P_{\rem{h_2}}} 
t_3 \wideoverunder{\rem{h_3}}{P_{\rem{h_3}}} 
t_4 \overunder{}{} 
\cdots \,. 
\end{equation} 
(cf.\ (\ref{eq:cyclic1}))  
The cluster at each point $t_m$ is given by  
$\textbf{x}(t_m)=\{y_h\,:\,h\in H_m\}$; 
as before, each cluster variable $y_h$ corresponds to the variable 
$x_{\rem{h}}$. 
The exchange polynomial $P_{\rem{h}}$ for an edge  
$~\bullet\overunder{\rem{h}}{} \bullet~$ with $h=(i,j)$ is given by 
\begin{equation}  
\label{eq:Ph-knight} 
P_{\rem{h}}=\alpha x_{\rem{(i,j-1)}}x_{\rem{(i-2,j)}} 
+\beta x_{\rem{(i-1,j)}} x_{\rem{(i-1,j-1)}}.  
\end{equation}  
Then equations (\ref{eq:elkies}) become 
the exchange relations associated with this pattern.  
 
To establish the Laurent phenomenon, we will  
complete the caterpillar pattern by attaching 
``legs'' to each vertex $t_m$  
and assigning exchange polynomials to these legs so that  
the appropriate analogues of conditions (\ref{eq:GEP1a}),  
(\ref{eq:GEP2a}) and (\ref{eq:GEP3a}) are satisfied.  
Since we now work over the polynomial ring  
$\AAA[x_a : a \in\mathcal{H}/\!\!\sim\,]$ 
in infinitely many indeterminates,  
the number of legs attached to every vertex $t_m$ will also be infinite 
(one for every label  
$a$ different from $\rem{h_{m-1}}$ and $\rem{h_{m}}$).  
This will not matter much for our argument though: to prove the Laurentness 
for any $y_{h_m}$, we will simply restrict our attention to the finite part  
of the infinite caterpillar tree lying between $t_0$ and  
$t_{\rm head} = t_{m+1}$, and to the legs labeled by $\rem{h_k}$  
for $0 \leq k \leq m$.

The role of conditions (\ref{eq:GEP1a}) and (\ref{eq:GEP2a}) 
is now played by the observation that each exchange  
polynomial $P_{\rem{h}}$ is not divisible by any variable~$x_a\,$,  
and furthermore every specialization 
$P_{\rem{h}}\bigl|_{x_a\leftarrow 0}$  
is an irreducible element of the Laurent polynomial ring. 
  
%We then complete the caterpillar pattern by attaching 
%(infinitely many) ``legs'' to each vertex $t_m$  
%and assigning exchange polynomials to these legs so that  
%the appropriate analogue of condition (\ref{eq:GEP3a}) is satisfied.  
%This will be done with the 
%help of a recursive procedure analogous to 
%(\ref{eq:G-step1})--(\ref{eq:G-step3}). 
 
To formulate the analogue of (\ref{eq:GEP3a}),  
let us fix an equivalence class $a\in\mathcal{H}/\!\!\sim\,$  
and concentrate on defining the exchange polynomials for the legs  
labeled by~$a$ and attached to the vertices squeezed between two  
consecutive occurrences of the label $a$ on the spine: 
 \vspace{-.5in} 
\begin{equation} 
\label{eq:knight-s} 
\bullet 
\overunder{a}{P_a} 
\bullet 
\overunder{a_1}{}\hspace{-.2in} 
\begin{array}{c}\\[.43in] \bullet\\ 
 \hspace{.12in}{\scriptstyle a} \Biggl| \scriptstyle G_1\ \\ 
\bullet\end{array} 
 \hspace{-.2in}\overunder{a_2}{} 
\hspace{-.2in} 
\begin{array}{c}\\[.43in] \bullet\\ \hspace{.08in}{\scriptstyle a} \Biggl|  
{\scriptstyle G_2}\\ 
\bullet\end{array} 
 \hspace{-.2in}\overunder{}{} 
\hspace{-.2in} 
\begin{array}{c}\\[.43in] \bullet\\ \hspace{.08in}{\scriptstyle a}  
\Biggl| {\hspace{.1in}\ }\\ 
\bullet\end{array} 
 \hspace{-.2in}\overunder{}{} 
\hspace{-.2in} 
\begin{array}{c}\\[.43in] \bullet\\ \hspace{.08in}{\scriptstyle a} \Biggl| {\hspace{.1in}\ }\\ 
\bullet\end{array} 
 \hspace{-.2in}\overunder{a_{N-2}}{} 
\hspace{-.35in} 
\begin{array}{c}\\[.43in] \bullet\\ \hspace{.22in}{\scriptstyle a} \Biggl|  
{\scriptstyle G_{N-2}}\\ 
\bullet\end{array} 
 \hspace{-.35in} 
\overunder{a_{N-1}}{} 
\bullet 
\overunder{a}{P_a} 
\bullet 
. 
\end{equation} 
We note that the labels $a_1,\dots,a_{N-1}\in\mathcal{H}/\!\!\sim\,$  
appearing on the spine between these two  
occurrences of $a$ are distinct. 
%It is convenient to set $a_0=a_N=a$. 
For $m=N-2,N-3,\dots,1$, we denote by  
$G_m$ the exchange polynomial to be associated with the $a$-labeled leg attached  
between the edges labeled $a_m$ and $a_{m+1}$  
(cf.~(\ref{eq:knight-s})).

The polynomials $G_m$ are defined with the 
help of a recursive procedure analogous to 
(\ref{eq:G-step1})--(\ref{eq:G-step3}).  
We initialize $G_{N-1}=P_a$, 
and obtain each $G_{m-1}$ from $G_m$, as follows. 
The step (\ref{eq:G-step1}) is replaced by  
\begin{equation*}  
%\label{eq:G-step1}  
\stackrel{\sim}{G}_{m-1} = G_m\bigl|_{x_{a_m}\leftarrow \frac{Q_m}{x_{a_m}}}  
\end{equation*}  
with  
\begin{equation}  
\label{eq:Qm-knight}  
Q_m=P_{a_m}\bigl|_{x_a\leftarrow 0} \,.  
\end{equation}  
We then compute $\stackrel{\approx}{G}_{m-1}$ and $G_{m-1}$  
exactly as in (\ref{eq:G-step2})--(\ref{eq:G-step3}). 
By the argument given in the proof of Theorem~\ref{th:laurent-cyclic},  
the equality $G_0=P_a$ would imply the desired Laurentness 
(cf.~(\ref{eq:GEP3a})). 
  
%[address the issue of infinite number of legs]  
  
To simplify %explicit  
computations, we denote the equivalence classes  
``surrounding'' $a$, as shown below: 
\begin{equation} 
\label{eq:elkies-table} 
\begin{array}{ccccccc} 
\cdots & \cdots & \cdots & \cdots & \cdots & \cdots & \cdots \\ 
\cdots & q & p & f & c & a & \cdots \\[.05in] 
\cdots & f & c & a & e & b & \cdots \\[.05in] 
\cdots & a & e & b & g & d& \cdots \\[.03in] 
\cdots & \cdots & \cdots & \cdots & \cdots & \cdots & \cdots  
\end{array} 
. 
\end{equation} 
In other words, if $a=\rem{(i,j)}$, then $b=\rem{(i,j-1)}$, 
$c=\rem{(i-1,j)}$, etc.  
With this notation, we can redraw the pattern (\ref{eq:knight-s}) as  
follows: 
\vspace{-.6in} 
\begin{equation} 
\label{eq:elkies-6steps} 
\bullet 
\overunder{a}{} 
\bullet 
\overunder{g}{}\cdots\hspace{-.25in} 
\begin{array}{c}\\[.43in] \bullet\\ \hspace{.2in}{}_{a} \Biggl| {}_{G_{k-1}}\\ 
\bullet\end{array} 
 \hspace{-.25in}\overunder{f}{} 
\hspace{-.2in} 
\begin{array}{c}\\[.43in] \bullet\\ \hspace{.08in}{}_{a} \Biggl| {}_{G_{k}}\\ 
\bullet\end{array} 
 \hspace{-.2in}\overunder{e}{} 
%\hspace{-.25in} 
\bullet 
% \hspace{-.25in} 
\overunder{d}{} 
%\hspace{-.2in} 
%\begin{array}{c}\\[.43in] \bullet\\ \hspace{.08in}{}_{a} \Biggl| {}_{G_{k+2}}\\ 
%\bullet\end{array} 
%\bullet  
% \hspace{-.2in}  
\cdots\hspace{-.25in} 
\begin{array}{c}\\[.43in] \bullet\\ \hspace{.2in}{}_{a} \Biggl| {}_{G_{\ell-1}}\\ 
\bullet\end{array} 
 \hspace{-.25in} 
\overunder{c}{} 
\hspace{-.2in} 
\begin{array}{c}\\[.43in] \bullet\\ \hspace{.08in}{}_{a} \Biggl| {}_{G_\ell}\\ 
\bullet\end{array} 
 \hspace{-.2in} 
\overunder{b}{} 
\bullet 
\cdots\overunder{a}{} 
\bullet 
,  
\end{equation} 
for appropriate values of $k$ and~$\ell$.

We will call a value of $m$ \emph{essential} if $G_{m-1}\neq G_m\,$.  
We are going to see that the essential values of $m$ 
%(i.e., those with $G_{m-1}\!\neq\! G_m$)  
are those for which $a_m\in\{b,c,e,f\}$;  
in the notation of (\ref{eq:elkies-6steps}), these values are  
$\ell+1$, $\ell$, $k+1$, and~$k$. 
%With this in mind, we compute (cf.\ (\ref{eq:Ph-knight})  
%and (\ref{eq:elkies-table}):  
% and (\ref{eq:Qm-knight})):  
%\begin{eqnarray} 
%\label{eq:elkies-a} 
%P_a &=& \alpha x_b x_f + \beta x_c x_e \,, \\ 
%\label{eq:elkies-b} 
%P_b &=& \alpha x_a x_d + \beta x_e x_g \,, \\ 
%\label{eq:elkies-c} 
%P_c &=& \alpha x_e x_p + \beta x_a x_f \,, \\ 
%\label{eq:elkies-e} 
%P_e &=& \alpha x_c x_g + \beta x_a x_b \,, \\ 
%\label{eq:elkies-f} 
%P_f &=& \alpha x_a x_q + \beta x_c x_p \,.  
%\end{eqnarray} 
  
We initialize $G_{N-1}=P_a=\alpha x_b x_f + \beta x_c x_e\,$.  
The values of $m$ in the interval $\ell<m<N$ are not essential since  
the variable $x_{a_m}$ does not enter $P_a$, which is furthermore  
not divisible by~$Q_m$ (because the latter involves variables absent  
in~$P_a$). 
  
The first essential value is $m=\ell+1$, with $a_m=b$:  
\[  
\begin{array}{l}  
Q_{\ell+1}=P_b|_{x_a\leftarrow 0}  
=(\alpha x_a x_d + \beta x_e x_g)|_{x_a\leftarrow 0}  
=\beta x_e x_g\,,\\[.1in]  
\stackrel{\sim}{G}_\ell 
= P_a\bigl|_{x_b\leftarrow \frac{Q_{\ell+1}}{x_b}} 
=\alpha \frac{\beta x_e x_g}{x_b} x_f + \beta x_c x_e\,,  
\\[.1in]  
{G}_\ell = \alpha x_g x_f + x_b x_c\,. 
\end{array}  
\] 
  
Step $m=\ell$ (here $a_m=c$):  
\[  
\begin{array}{l}  
Q_\ell=P_c|_{x_a\leftarrow 0}=(\alpha x_e x_p + \beta x_a x_f)|_{x_a\leftarrow 0}=  
\alpha x_e x_p \,,\\[.1in]  
\stackrel{\sim}{G}_{\ell-1}  
=G_\ell \bigl|_{x_c\leftarrow \frac{Q_\ell}{x_c}} 
=\alpha x_g x_f + x_b \frac{\alpha x_e x_p}{x_c} \,,\\[.1in]  
{G}_{\ell-1}=x_c x_g x_f + x_b x_e x_p \,. 
\end{array}  
\] 
Notice that ${G}_{\ell-1}$ does not involve $x_d$, so 
the value $m=k+2$ is not essential, as are the rest of the values in the interval  
$k+1<m<\ell$. 
  
Step $m=k+1$, with $a_m=e$:  
\[  
\begin{array}{l}  
Q_{k+1}=P_e|_{x_a\leftarrow 0}  
=(\alpha x_c x_g + \beta x_a x_b)|_{x_a\leftarrow 0}=\alpha x_c  
x_g\,,\\[.1in]  
\stackrel{\sim}{G}_k  
%=G_{k+1}\bigl|_{x_e\leftarrow \frac{Q_0}{x_e}} 
=x_c x_g x_f + x_b x_p \frac{\alpha x_c x_g}{x_e}\,, \\[.1in]  
{G}_k=x_f x_e + \alpha x_b x_p\,. 
\end{array}  
\] 
  
Step $m=k$, with $a_m=f$:  
\[  
\begin{array}{l}  
Q_k=P_f|_{x_a\leftarrow 0}=(\alpha x_a x_q + \beta x_c  
x_p)|_{x_a\leftarrow 0}=\beta x_c x_p \,,\\[.1in]  
\stackrel{\sim}{G}_{k-1}  
=\frac{\beta x_c x_p}{x_f} x_e + \alpha x_b x_p\,, \\[.1in]  
{G}_{k-1}=\beta x_c x_e + \alpha x_b x_f\,. 
\end{array}  
\] 
The values of $m$ in the interval $0<m<k$ are not essential  
since none of the corresponding variables $x_{a_m}$ appears in 
${G}_{k-1}$; in particular, $m=1$ is not essential,  
since ${G}_{k-1}$ does not involve~$x_g\,$. 
Hence 
\[  
G_0={G}_{k-1}=\beta x_c x_e + \alpha x_b x_f=P_a\,,  
\]  
as desired. The Laurentness is proved.

\begin{remark} 
\label{remark:initial-conditions}
{\rm  
The Laurent phenomenon for the recurrence (\ref{eq:elkies}) 
actually holds in greater generality. 
Specifically, one can replace $\mathcal{H}$ by any subset 
of $\ZZ^2$  
which satisfies the following analogues of conditions
(\ref{eq:cube-h1})--(\ref{eq:cube-h2}) and (\ref{eq:oct-h1})--(\ref{eq:oct-h2}): 
\begin{eqnarray} 
\label{eq:h-axiom1} 
&&\text{if $h\in\mathcal{H}$, then $h'\in\mathcal{H}$ 
whenever $h\leq h'$;} 
\\ 
\label{eq:h-axiom2} 
&&\text{for any $h'\in\mathcal{H}$, the set 
$\{h\in\mathcal{H}:h\leq h'\}$ is finite.} \hspace{1.3in}
\end{eqnarray}  
Then take %$H_{\rm init}$ by the set  
$
H_{\rm init} =\{h\in\mathcal{H}\,:\,h^-\notin\mathcal{H}\}$.  

The proof of Laurentness only needs one adjustment,
concerning the choice of a linear extension~$\prec\,$.
Specifically, while proving that $y_h$ is given by a Laurent polynomial,
take a finite set $\mathcal{H}^{(h)}\subset \mathcal{H}$ 
containing $h$ and satisfying the conditions
\begin{eqnarray} 
&&\text{if $h'\in\mathcal{H}^{(h)}$, then $h''\in\mathcal{H}^{(h)}$ 
whenever $h''\leq h'$ and $h''\in\mathcal{H}$;} 
\\ 
&&\text{for any $h'\in\mathcal{H}$ such that $h'\leq h$, there
exists $h''\in\mathcal{H}^{(h)}$ such that}\\
\nonumber
&&\text{$h''\geq h$ and $h''\sim h$.}
\end{eqnarray} 
(The existence of $\mathcal{H}^{(h)}$ follows from
(\ref{eq:h-axiom1})--(\ref{eq:h-axiom2}).) 
Then define $\preceq$ exactly as before on the set 
$\mathcal{H}^{(h)}$;
set $h'\prec h''$ for any $h'\in \mathcal{H}^{(h)}$
and $h''\in \mathcal{H}-\mathcal{H}^{(h)}$;
and define $\preceq$ on the complement $\mathcal{H}-\mathcal{H}^{(h)}$ 
by an arbitrary linear extension of~$\leq$. 
These conditions ensure that the sets $H_m$ needed in the proof of
Laurentness of the given $y_h$ are well defined,
and that the rest of the proof proceeds smoothly. 
%The rest of the examples in this section can also be generalized in this 
%fashion.  
} 
\end{remark} 

Armed with the techniques developed above in this section,
we will now prove the main theorems stated in the introduction. 

\medskip
 
\noindent
\textbf{Proof of Theorem~\ref{th:kleber}.}
Our argument is parallel to that in Example~\ref{example:elkies}, 
so we skip the steps which are identical in both proofs.  
For simplicity of exposition, we present the proof in the special case
$\mathcal{H}=\ZZ_{\geq 0}^3$;
the case of general $\mathcal{H}$ requires the same adjustments
as those described in Remark~\ref{remark:initial-conditions}. 
 
We define the product partial order $\leq$ and a compatible linear 
order $\preceq$ on $\mathcal{H}$ by  
\begin{eqnarray*} 
(i_1,j_1,k_1)\leq (i_2,j_2,k_2)  
&\!\stackrel{\rm def}{\Leftrightarrow}\!& 
(i_1\leq i_2)\text{~and~} (j_1\leq j_2)\text{~and~} (k_1\leq k_2),\\ 
(i_1,j_1,k_1)\preceq (i_2,j_2,k_2) &\!\stackrel{\rm def}{\Leftrightarrow}\!& 
(i_1+j_1+k_1<i_2+j_2+k_2)\\ 
&&\text{~or~}  
(i_1+j_1+k_1= i_2+j_2+k_2 \text{~and~} i_1+j_1< i_2+j_2)\\ 
&&\text{~or~}  
(i_1+j_1= i_2+j_2 \text{~and~} k_1=k_2 \text{~and~} i_1\leq i_2).  
\end{eqnarray*} 
For $h=(i,j,k)$, we set $h^-=(i-1,j-1,k-1)$; thus,  
the exchange relation (\ref{eq:cube-frac})  
expresses the product $y_h\cdot y_{h^-}$ as a polynomial in the 
variables $y_{h'}\,$, for $h^-< h'< h$. 
 
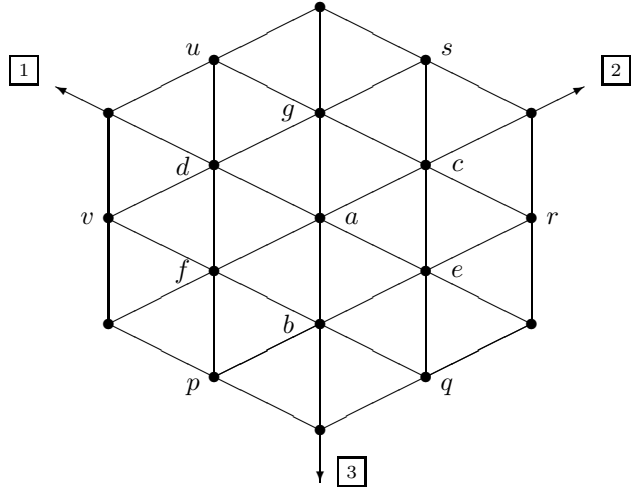
\begin{figure}[ht] 
\setlength{\unitlength}{4pt} 
 
\centering 
\begin{picture}(40,45)(0,-4) 
%\thicklines 
%\setlength{\linewidth}{3pt} 
\multiput(0,30)(10,5){3}{\circle*{1}} 
\multiput(0,20)(10,5){4}{\circle*{1}} 
\multiput(0,10)(10,5){5}{\circle*{1}} 
\multiput(10,5)(10,5){4}{\circle*{1}} 
\multiput(20,0)(10,5){3}{\circle*{1}} 
 
\put(0,30){\line(2,1){20}} 
\put(0,20){\line(2,1){30}} 
\put(0,10){\line(2,1){40}} 
\put(10,5){\line(2,1){30}} 
\put(20,0){\line(2,1){20}} 
\put(10,5){\line(2,1){10}} 
\put(30,5){\line(2,1){10}} 
 
\put(0,10){\line(0,1){20}} 
\put(10,5){\line(0,1){30}} 
\put(20,0){\line(0,1){40}} 
\put(30,5){\line(0,1){30}} 
\put(40,10){\line(0,1){20}} 
 
\put(0,10){\line(2,-1){20}} 
\put(0,20){\line(2,-1){30}} 
\put(0,30){\line(2,-1){40}} 
\put(10,35){\line(2,-1){30}} 
\put(20,40){\line(2,-1){20}} 
 
\put(40,30){\vector(2,1){5}} 
\put(0,30){\vector(-2,1){5}} 
\put(20,0){\vector(0,-1){5}} 
 
\put(-2,29){\makebox(0,0){}} 
\put(-2,20){\makebox(0,0){$v$}} 
\put(-2,10){\makebox(0,0){}} 
\put(8,36){\makebox(0,0){$u$}} 
\put(7,25){\makebox(0,0){$d$}} 
\put(7,15){\makebox(0,0){$f$}} 
\put(8,4){\makebox(0,0){$p$}} 
\put(20,42){\makebox(0,0){}} 
\put(17,30){\makebox(0,0){$g$}} 
\put(23,20){\makebox(0,0){$a$}} 
\put(17,10){\makebox(0,0){$b$}} 
\put(18,-1){\makebox(0,0){}} 
\put(32,36){\makebox(0,0){$s$}} 
\put(33,25){\makebox(0,0){$c$}} 
\put(33,15){\makebox(0,0){$e$}} 
\put(32,4){\makebox(0,0){$q$}} 
\put(42,29){\makebox(0,0){}} 
\put(42,20){\makebox(0,0){$r$}} 
\put(42,10){\makebox(0,0){}} 
 
\put(-8,34){\makebox(0,0){$\boxed{\scriptstyle 1}$}} 
\put(48,34){\makebox(0,0){$\boxed{\scriptstyle 2}$}} 
\put(23,-4){\makebox(0,0){$\boxed{\scriptstyle 3}$}} 
 
\end{picture} 
 
\caption{The cube recurrence} 
\label{fig:kleber-classes} 
\end{figure}

All the steps in Example~\ref{example:elkies} leading to the creation of the  
infinite exchange pattern (\ref{eq:knight-spine})  
are repeated verbatim.  
Instead of (\ref{eq:Ph-knight}), the exchange polynomials  
$P_{\rem{h}}$ along the spine are now given by 
\[ 
\begin{array}{l} 
\!\!P_{\rem{(i,j,k)}}\\[.1in] 
\!\!=\alpha x_{\rem{(i-1,j,k)}}x_{\rem{(i,j-1,k-1)}} 
\!+\!\beta x_{\rem{(i,j-1,k)}}x_{\rem{(i-1,j,k-1)}} 
\!+\!\gamma x_{\rem{(i,j,k-1)}}x_{\rem{(i-1,j-1,k)}}. 
\end{array} 
\] 
%Referring to Figure~\ref{fig:kleber-classes} and assuming 
%$[(i,j,k)]=r$,  
%this would give, for example, 
%\[ 
%P_r=x_g x_s + x_m x_q + x_v x_l\,.  
%\] 
%The rest of the argument follows the footsteps of Example~\ref{example:elkies}.  
The role of (\ref{eq:elkies-table}) is now played by
Figure~\ref{fig:kleber-classes}, 
which shows the ``vicinity'' of an equivalence class~$a$.
This figure displays the orthogonal projection of $\mathcal{H}$ along 
the vector $(1,1,1)$. 
Thus the vertices represent equivalence classes in 
$\mathcal{H}/\!\!\sim\,$. 
For example, if $a=\rem{(i,j,k)}$, then 
\begin{equation*} 
\begin{array}{lll} 
b=\rem{(i,j,k-1)}, &c=\rem{(i,j-1,k)}, &d=\rem{(i-1,j,k)},\\[.1in] 
e=\rem{(i,j-1,k-1)}, &f=\rem{(i-1,j,k-1)}, &g=\rem{(i-1,j-1,k)}.  
\end{array} 
\end{equation*} 
With this notation, we have:  
\[ 
P_a=\alpha x_d x_e + \beta x_c x_f + \gamma x_b x_g\,.  
\] 
 
With the polynomials $G_1,G_2,\dots$ defined as in 
(\ref{eq:knight-s}),  
the essential values of $m$ are now those for which 
$a_m\in\{b,c,d,e,f,g\}$.  
(The verification that the rest of the values are not essential is 
left to the reader.)  
We denote these values by $m_1,\dots,m_6\,$, respectively.  
 
The computation of the polynomials $G_m$ begins by initializing 
\[ 
G_{N-1}=P_a=\alpha x_d x_e + \beta x_c x_f + \gamma x_b x_g\,. 
\] 
Step $m=m_1$, $a_m=b$: 
\[ 
\begin{array}{rl} 
Q_{m_1}=&\!\!\!P_{b}|_{x_a\leftarrow 0} 
=\alpha x_f x_q + \beta x_e x_p\,; 
\\[.09in] 
\stackrel{\sim}{G}_{m_1-1}=&\!\!\!G_{m_1}\bigl|_{x_{b}\leftarrow 
 \frac{Q_{m_1}}{b}} 
=\alpha x_d x_e + \beta x_c x_f + \gamma \frac{\alpha x_f x_q + \beta x_e x_p}{x_b} x_g 
 \,; 
\\[.09in] 
G_{m_1-1}=&\!\!\! 
\alpha x_b x_d x_e + \beta x_b x_c x_f  
+ \alpha \gamma x_f x_g x_q + \beta \gamma x_e x_g x_p \,. 
\end{array} 
\] 
Step $m=m_2$, $a_m=c$: 
\[ 
\begin{array}{rl} 
Q_{m_2}=&\!\!\!P_c|_{x_a\leftarrow 0} 
= \alpha x_g x_r + \gamma x_e x_s;\\[.09in] 
\stackrel{\sim}{G}_{m_2-1}=&\!\!\! 
\alpha x_b x_d x_e + \beta x_b \frac{\alpha x_g x_r + \gamma x_e x_s}{x_c} x_f  
+ \alpha \gamma x_f x_g x_q + \beta \gamma x_e x_g x_p 
\,;  
\\[.09in] 
G_{m_2-1}=&\!\!\! 
\alpha x_b x_c x_d x_e  
\!+\! \alpha \beta x_b x_f x_g x_r \!+\! \beta \gamma x_b x_e x_f x_s 
\!+\! \alpha \gamma x_c x_f x_g x_q \!+\! \beta \gamma x_c x_e x_g x_p 
\,.  
\end{array} 
\] 
Step $m=m_3$, $a_m=d$: 
\[ 
\begin{array}{rl} 
Q_{m_3}=&\!\!\!P_d|_{x_a\leftarrow 0} 
= \beta x_g x_v + \gamma x_f x_u\,; 
\\[.09in] 
\stackrel{\sim}{G}_{m_3-1}=&\!\!\! 
\alpha x_b x_c \frac{\beta x_g x_v + \gamma x_f x_u}{x_d} x_e \\[.09in] 
& 
+ \alpha \beta x_b x_f x_g x_r + \beta \gamma x_b x_e x_f x_s 
+ \alpha \gamma x_c x_f x_g x_q + \beta \gamma x_c x_e x_g x_p 
\,; 
\\[.09in] 
G_{m_3-1}=&\!\!\! 
 \alpha \beta x_b x_c x_e x_g x_v + \alpha \gamma x_b x_c x_e x_f x_u 
+ \beta \gamma x_b x_d x_e x_f x_s + \beta \gamma x_c x_d x_e x_g x_p \\[.09in] 
& 
+ \alpha \beta x_b x_d x_f x_g x_r + \alpha \gamma x_c x_d x_f x_g x_q  
\,. 
\end{array} 
\] 
Step $m=m_4$, $a_m=e$: 
\[ 
\begin{array}{rl} 
Q_{m_4}=&\!\!\!P_e|_{x_a\leftarrow 0} 
= \beta x_b x_r + \gamma x_c x_q\,; 
\\[.09in] 
\stackrel{\sim}{G}_{m_4-1}=&\!\!\! 
 \frac{Q_{m_4}}{x_e}(\alpha \beta x_b x_c x_g 
 x_v + \alpha \gamma x_b x_c x_f x_u + \beta \gamma x_b x_d x_f x_s + 
\beta \gamma x_c x_d x_g x_p) \\[.09in]  
& 
+ \alpha x_d x_f x_g Q_{m_4} 
\,;  
\\[.09in] 
G_{m_4-1}=&\!\!\! 
\alpha \gamma x_b x_c x_f x_u + \beta \gamma x_b x_d x_f x_s  
+ \alpha x_d x_e x_f x_g \\[.09in] 
& 
+\alpha \beta x_b x_c x_g x_v + \beta \gamma x_c x_d x_g x_p  
\,. 
\end{array} 
\] 
Step $m=m_5$, $a_m=f$: 
\[ 
\begin{array}{rl} 
Q_{m_5}=&\!\!\!P_f|_{x_a\leftarrow 0} 
= \alpha x_b x_v + \gamma x_d x_p \,; 
\\[.09in] 
\stackrel{\sim}{G}_{m_5-1}=&\!\!\! 
\frac{Q_{m_5}}{x_f} 
(\alpha \gamma x_b x_c x_u + \beta \gamma x_b x_d x_s  
+ \alpha x_d x_e x_g) +\beta x_c x_g Q_{m_5} 
\,; 
\\[.09in] 
G_{m_5-1}=&\!\!\! 
 \alpha x_d x_e x_g +\beta x_c x_f x_g  
+\alpha \gamma x_b x_c x_u + \beta \gamma x_b x_d x_s  
\,. 
\end{array} 
\] 
Step $m=m_6$, $a_m=g$: 
\[ 
\begin{array}{rl} 
Q_{m_6}=&\!\!\!P_g|_{x_a\leftarrow 0} 
= \alpha x_c x_u + \beta x_d x_s \,;  
\\[.09in] 
\stackrel{\sim}{G}_{m_6-1}=&\!\!\! 
\frac{Q_{m_6}}{x_g}(\alpha x_d x_e +\beta x_c x_f)  
+\gamma x_b Q_{m_6} 
\,; 
\\[.09in] 
G_{m_6-1}=&\!\!\! 
\alpha x_d x_e +\beta x_c x_f +\gamma x_b x_g  
=P_a 
\,,  
\end{array} 
\] 
%proving Theorem~\ref{th:kleber}.  
completing the proof. 
\endproof
 
We will now deduce the Gale-Robinson conjecture from 
Theorem~\ref{th:kleber}. 
 
\medskip 
 
\noindent 
\textbf{Proof of Theorem~\ref{th:gale-robinson}.} 
To prove the Laurentness of a given element $y_N$ of the Gale-Robinson
sequence $(y_m)$, we define the array 
$(z_{ijk})_{(i,j,k)\in\mathcal{H}}$ by setting
$z_{ijk}=y_{N+pi+qj+rk}\,$,
with the indexing set 
\[
\mathcal{H}=\mathcal{H}(N)=\{(i,j,k)\in\ZZ^3\,:\,N+pi+qj+rk\geq 0
\}\,.
\]
Then (\ref{eq:gale-robinson}) implies 
that the $z_{ijk}$ satisfy the cube recurrence (\ref{eq:cube-frac}).
Note that $\mathcal{H}$ satisfies the conditions 
(\ref{eq:cube-h1})--(\ref{eq:cube-h2}).
Thus Theorem~\ref{th:kleber} applies to $(z_{ijk})$, with  
%and the corresponding initial set $H_{\rm init}$ is given by
$H_{\rm init}=\{(a,b,c)\in\ZZ^3: 0\leq N+pa+qb+rc<n \}$.
It remains to note that $y_N=z_{000}\,$, 
while for any $(a,b,c)\in H_{\rm init}\,$,
we have $z_{abc}=y_m$ with $0\leq m<n$.
\endproof

\medskip 
 
\noindent 
\textbf{Proof of Theorem~\ref{th:oct}.} 
This theorem is proved by the same argument as
Theorem~\ref{th:kleber}. 
We treat the Mills-Robbins-Rumsey special case
(\ref{eq:half-lattice-mod2}) (cf.\ also (\ref{eq:lattice-mod2})); 
similarly to Theorem~\ref{th:kleber}, 
the case of general $\mathcal{H}$ requires the standard adjustments
described in Remark~\ref{remark:initial-conditions}. 
We use the partial order on the lattice $L$ defined by
\[
(i,j,k)\leq (i',j',k') \,:\, 
|i'-i|+|j'-j|\leq k'-k\,.
\] 
For $h=(i,j,k)\in L$, we set $h^-=(i,j,k-2)$,
and define the equivalence relation~$\sim$ accordingly. 
Figure~\ref{fig:octahedral} shows equivalence classes ``surrounding'' a given 
class~$a$ (cf.\ Figure~\ref{fig:kleber-classes}).

\begin{figure}[ht]  
\setlength{\unitlength}{6pt}  
  
\centering  
\begin{picture}(20,11)(0,0)  
%\thicklines  
%\setlength{\linewidth}{3pt}  
  
\put(5,0){\circle*{0.7}}  
\put(10,0){\circle{0.7}}  
\put(15,0){\circle*{0.7}}  
\put(5,5){\circle{0.7}}  
\put(10,5){\circle*{0.7}}  
\put(15,5){\circle{0.7}}  
\put(5,10){\circle*{0.7}}  
\put(10,10){\circle{0.7}}  
\put(15,10){\circle*{0.7}}  
  
\put(5,0){\line(1,0){10}}  
\put(5,5){\line(1,0){10}}  
\put(5,10){\line(1,0){10}}  
\put(5,0){\line(0,1){10}}  
\put(10,0){\line(0,1){10}}  
\put(15,0){\line(0,1){10}}  
  
\put(0,5){\vector(0,1){3}}  
\put(20,3){\vector(1,0){3}}  
  
\put(6,6){\makebox(0,0){$d$}}  
\put(11,6){\makebox(0,0){$a$}}  
\put(16,6){\makebox(0,0){$c$}}  
\put(11,1){\makebox(0,0){$e$}}  
\put(6,11){\makebox(0,0){$p$}}  
\put(11,11){\makebox(0,0){$b$}}  
\put(16,11){\makebox(0,0){$q$}}  
\put(6,1){\makebox(0,0){$s$}}  
\put(16,1){\makebox(0,0){$r$}}  
  
\put(-1,7){\makebox(0,0){$j$}}  
\put(22,4){\makebox(0,0){$i$}}  
  
\end{picture}  
  
\caption{}  
\label{fig:octahedral}  
\end{figure}  
  
The initialization polynomial $G_{N-1}=P_a$
is given by  $P_a=\alpha x_c x_d + \beta x_b x_e\,$. 
The table below displays $a_m$, $Q_m$, $\stackrel{\sim}{G}_{m-1}$, 
and $G_{m-1}$ for all essential values of~$m$.  

\[  
\begin{array}{cccc}  
a_m & Q_m &\stackrel{\sim}{G}_{m-1} & G_{m-1} \\[.1in] 
\hline\\ 
b & \alpha x_p x_q & \alpha x_c x_d +\alpha \beta \frac{x_p x_q}{x_b} x_e  
& x_b x_c x_d + \beta x_e x_p x_q 
\\[.1in]  
c & \beta x_q x_r & \beta \frac{x_q x_r}{x_c} x_b x_d + \beta x_e x_p x_q 
& x_b x_d x_r + x_c x_e x_p  
\\[.1in]  
d & \beta x_p x_s & \beta \frac{x_p x_s}{x_d} x_b x_r + x_c x_e x_p 
& \beta x_b x_r x_s + x_c x_d x_e  
\\[.1in]  
e & \alpha x_r x_s & \beta x_b x_r x_s + \alpha \frac{x_r x_s}{x_e} x_c x_d  
& \beta x_b x_e + \alpha x_c x_d %= P_a
\end{array}  
\]  
We see that $G_0=G_{e-1}=P_a\,$, completing the proof.
\endproof

\medskip 
 
\noindent 
\textbf{Proof of Theorem~\ref{th:two-term-gale-robinson}.} 
The proof mimics the above proof of Theorem~\ref{th:gale-robinson}.
To prove the Laurentness of an element $y_N$ of the 
sequence $(y_m)$ satisfying (\ref{eq:gale-robinson-two-term}), we
define the array 
$(z_{ijk})_{(i,j,k)\in\mathcal{H}}$ by setting
$z_{ijk}=y_{N+\ell(i,j,k)}\,$,
where $\ell(i,j,k)=n\frac{i+j+k}{2}-pi-qj$. 
The indexing set $\mathcal{H}$ is now given by 
\[
\mathcal{H}=\mathcal{H}(N)=\{(i,j,k)\in\ZZ^3\,:\,N+\ell(i,j,k)\geq 0
\}\,.
\]
Then (\ref{eq:gale-robinson-two-term}) implies 
that the $z_{ijk}$ satisfy the octahedron recurrence (\ref{eq:oct}).
It is easy to check that $\mathcal{H}$ satisfies the conditions 
(\ref{eq:oct-h1})--(\ref{eq:oct-h2}).
Thus Theorem~\ref{th:oct} applies to $(z_{ijk})$, with  
%and the corresponding initial set $H_{\rm init}$ is given by
$H_{\rm init}=\{(a,b,c)\in L\,:\, 0\leq N+\ell(a,b,c)<n \}$,
and the theorem follows.
\endproof

We conclude this section by a couple of examples
in which the Laurent phenomenon is established by the same technique
as above.  
In each case, we provide: 
\begin{itemize} 
\item  
a picture of the equivalence classes ``surrounding'' a given 
class~$a$, 
which plays the role of (\ref{eq:elkies-table}) in 
Example~\ref{example:elkies};  
\item 
the initialization polynomial $G_{N-1}=P_a$; 
\item 
a table showing $a_m$, $Q_m$, $\stackrel{\sim}{G}_{m-1}$, 
and $G_{m-1}$ for all essential values of~$m$.  
\end{itemize} 
 
\begin{example}  
\label{example:frieze-patterns}  
 (Frieze patterns)  
{\rm 
The generalized frieze pattern recurrence 
(cf., e.g., \cite{conway-coxeter, ec2}) is  
\begin{equation}  
y_{ij} y_{i-1,j-1}  
= \varepsilon\, y_{i,j-1} y_{i-1,j} + \beta\,, 
\end{equation}  
where $\varepsilon\in\{1,-1\}$.  
To prove Laurentness (over $\ZZ[\beta]$),  
refer to Figure~\ref{fig:frieze}. 
Then $P_a=\varepsilon\,x_b\, x_c +\beta$, and the essential steps are:  
\[  
\begin{array}{cccc}  
a_m & Q_m &\stackrel{\sim}{G}_{m-1} & G_{m-1} \\[.1in] 
\hline\\ 
b & \beta & \frac{\varepsilon\,\beta\, x_c}{x_b} + \beta  
&  
\varepsilon\,x_c + x_b  
\\[.1in]  
c & \beta  
&  
\frac{\varepsilon\,\beta}{x_c}+ x_b  
&  
\beta+\varepsilon^{-1} x_b \,x_c %=P_a\,. 
\end{array}  
\]  
  
\begin{figure}[ht]  
\setlength{\unitlength}{4pt}  
  
\centering  
\begin{picture}(20,12)(0,0)  
%\thicklines  
%\setlength{\linewidth}{3pt}  
  
\put(10,0){\circle*{1}}  
\put(5,5){\circle*{1}}  
\put(15,5){\circle*{1}}  
\put(10,10){\circle*{1}}  
  
\put(5,5){\line(1,1){5}}  
\put(15,5){\line(-1,1){5}}  
  
\put(10,0){\vector(-1,1){10}}  
\put(10,0){\vector(1,1){10}}  
  
\put(7,5){\makebox(0,0){$c$}}  
\put(17,5){\makebox(0,0){$b$}}  
\put(12,0){\makebox(0,0){$a$}}  
\put(12,10){\makebox(0,0){$a$}}  
  
\put(-1,9){\makebox(0,0){$i$}}  
\put(21,9){\makebox(0,0){$j$}}  
  
\end{picture}  
  
\caption{}  
\label{fig:frieze}  
\end{figure}  
}  
\end{example}  
  
\begin{example}  
\label{example:number-walls}  
 (Number walls)  
{\rm  
Consider the 2-dimensional recurrence  
\begin{equation}  
\label{eq:number-walls} 
y_{ij}y_{i,j-2}  
=y_{i-1,j-1}^p y_{i+1,j-1}^r + y_{i,j-1}^q\,,  
\end{equation}  
where $p$, $q$, and $r$ are nonnegative integers. 
To prove Laurentness,  
refer to Figure~\ref{fig:number-wall}. 
Then $P_a=x_d^p x_b^r + x_c^q$, 
and the essential steps are:  
\[  
\begin{array}{cccc}  
a_m & Q_m &\stackrel{\sim}{G}_{m-1} & G_{m-1} \\[.1in] 
\hline\\ 
b &  
x_f^q  
&  
x_d^p \bigl(\frac{x_f^q}{x_b}\bigr)^r + x_c^q  
&  
x_d^p x_f^{qr} +x_c^q x_b^r  
\\[.1in]  
c & x_g^p x_f^r  
&  
x_d^p x_f^{qr} 
+\bigl(\frac{x_g^p x_f^r}{x_c}\bigr)^q x_b^r  
&  
x_d^p x_c^q +x_g^{pq} x_b^r  
\\[.1in]  
d & x_g^q  
&  
\bigl(\frac{x_g^q}{x_d}\bigr)^p x_c^q +x_g^{pq} x_b^r  
&  
x_c^q +x_b^r x_d^p  
\end{array}  
\]  
  
\begin{figure}[ht]  
\setlength{\unitlength}{6pt}  
  
\centering  
\begin{picture}(20,11)(0,0)  
%\thicklines  
%\setlength{\linewidth}{3pt}  
  
\put(5,0){\circle*{0.7}}  
\put(10,0){\circle*{0.7}}  
\put(15,0){\circle*{0.7}}  
\put(5,5){\circle*{0.7}}  
\put(10,5){\circle*{0.7}}  
\put(15,5){\circle*{0.7}}  
\put(10,10){\circle*{0.7}}  
  
\put(5,0){\line(1,0){10}}  
\put(5,5){\line(1,0){10}}  
\put(5,0){\line(0,1){5}}  
\put(10,0){\line(0,1){10}}  
\put(15,0){\line(0,1){5}}  
  
\put(0,5){\vector(0,1){3}}  
\put(20,3){\vector(1,0){3}}  
  
\put(6,6){\makebox(0,0){$d$}}  
\put(11,6){\makebox(0,0){$c$}}  
\put(16,6){\makebox(0,0){$b$}}  
\put(11,1){\makebox(0,0){$a$}}  
\put(11,10){\makebox(0,0){$a$}}  
\put(6,1){\makebox(0,0){$g$}}  
\put(16,1){\makebox(0,0){$f$}}  
  
\put(-1,7){\makebox(0,0){$j$}}  
\put(22,4){\makebox(0,0){$i$}}  
  
\end{picture}  
  
\caption{}  
\label{fig:number-wall}  
\end{figure}  
}  
\end{example}  
 
\begin{remark} 
{\rm 
As pointed out by J.~Propp, the Laurent phenomenon for certain 
special cases of Examples~\ref{example:frieze-patterns} 
and~\ref{example:number-walls} can be obtained by specialization 
of Example~\ref{example:octahedral}. 
} 
\end{remark}

\section{Homogeneous exchange patterns} 
\label{sec:homogeneous-exchange-patterns} 
 
In this section, we deduce Theorem~\ref{th:conway-scott} and a number of 
similar results from the following corollary of Theorem~\ref{th:laurent-gcd}. 
 
\begin{corollary} 
\label{cor:homogeneous-pattern} 
Let $\AAA$ be a unique factorization domain.  
Assume that a collection of nonzero polynomials $P_1,\dots,P_n\in\AAA[x_1,\dots,x_n]$ 
satisfies the following conditions: 
\begin{eqnarray} 
\label{eq:hom1} 
&&\text{Each $P_k$ does not depend on $x_k$, and is not divisible by any 
$x_i$, $i\in [n]$.}\\ 
\label{eq:hom2} 
&&\text{For any $i\neq j$, 
the polynomials $P_{ji}\!\stackrel{\rm def}{=}\!(P_j)|_{x_i=0}$ 
and $P_i$ are coprime.}\\ 
\label{eq:hom3} 
&&\text{For any $i\neq j$, we have}\\ 
\nonumber 
%\label{eq:P=LR0/Q0} 
&& \hspace{1in} L \cdot P_{ji}^b \cdot P_i\!=\! 
P_i\bigl|_{x_j\leftarrow \frac{P_{ji}}{x_j}}\,,\\ 
\nonumber 
&&\text{where $b$ is a nonnegative integer, and $L$ is a Laurent monomial whose }\\ 
\nonumber 
&& \text{coefficient lies in $\AAA$ and 
is coprime with~$P_i$.}  
\end{eqnarray} 
Let us define the rational transformations $F_i$, $i\in [n]$, by  
\begin{equation*} 
F_i : (x_1,\dots,x_n) \mapsto  
(x_1,\dots,x_{i-1},\displaystyle\frac{P_i}{x_i},x_{i+1},\dots,x_n). 
\end{equation*} 
Then any composition of the form $F_{i_1}\circ\cdots\circ F_{i_m}$  
is given by Laurent polynomials with coefficients in $\AAA$.  
\end{corollary} 
 
\proof 
Let $\TT_n$ denote a regular tree of degree $n$  
whose edges are labeled by elements of $[n]$ so that all edges incident to a 
given vertex have different labels. 
Assigning $P_i$ as an exchange polynomial for every edge of $\TT_n$ labeled by $i$, 
we obtain a ``homogeneous'' exchange pattern on $\TT_n$ satisfying 
conditions (\ref{eq:GEP1})--(\ref{eq:GEP3}) in Theorem~\ref{th:laurent-gcd}. 
This implies the desired Laurentness.  
%Since every path of size $m$ in $\TT_n$ can be included as the spine into 
%a finite caterpillar tree $\TT_{n,m} 
\endproof 
 
\begin{example} 
{\rm 
Let $n \geq 3$ be an integer, 
and let $P$ be a quadratic form given by 
$$P(x_1, \dots, x_n) = x_1^2 + \cdots + x_n^2 + 
\displaystyle \sum_{i < j} \alpha_{ij} x_i x_j \ .$$ 
Theorem~\ref{th:conway-scott} is a special case of 
Corollary~\ref{cor:homogeneous-pattern} 
for $P_i = P\bigl|_{x_i = 0}$ and $\AAA = \ZZ[\alpha_{ij}:i<j]$. 
Conditions (\ref{eq:hom1})--(\ref{eq:hom2}) are clear. 
To verify (\ref{eq:hom3}), note that 
$$P_i = P_{ji} + x_j^2 + x_j \displaystyle\Bigl(\sum_k \alpha_{kj} x_k 
+ \sum_\ell \alpha_{j\ell} x_l\Bigr) \ ,$$ 
where $k$ (resp.~$\ell$) runs over all indices such that 
$k \neq i$ and $k < j$ (resp.~$\ell \neq i$ and $\ell > j$). 
It follows that 
$$P_i\bigl|_{x_j\leftarrow \frac{P_{ji}}{x_j}} = 
P_{ji} + \frac{P_{ji}^2}{x_j^2} + \frac{P_{ji}}{x_j} 
\displaystyle\Bigl(\sum_k \alpha_{kj} x_k + \sum_\ell \alpha_{j\ell} x_l\Bigr) 
= \frac{P_{ji}}{x_j^2} P_i \ ,$$ 
verifying (\ref{eq:hom3}). 
} 
\end{example}

In the remainder of this section, we list a few more applications 
of Corollary~\ref{cor:homogeneous-pattern}. 
In each case, the verification of its conditions  
is straightforward.  
 
\begin{example} 
{\rm 
Let $P$ and $Q$ be monic palindromic 
%(i.e., self-reciprocal) 
polynomials in one variable: 
\[ 
\begin{array}{l} 
P(x)=(1+x^d)+\alpha_1(x+x^{d-1})+\alpha_2(x^2+x^{d-2})+\dots\,;\\[.1in] 
Q(x)=(1+x^e)+\beta_1 (x+x^{e-1})+\beta_2 (x^2+x^{e-2})+\dots\,. 
\end{array} 
\] 
% Then the polynomials 
% \[ 
% \begin{array}{l} 
% P(x_2)=\mu^2 \tilde P(\frac{x_2}{\lambda})\\ 
% Q(x_1)=\lambda^2 \tilde Q(\frac{x_1}{\mu}) 
% \end{array} 
% \] 
% define a homogeneous generalized exchange pattern of rank~2. 
Then every member of the sequence $y_0,y_1,\dots$ defined by the recurrence 
\[ 
y_k = 
\begin{cases} 
\displaystyle\frac{\mu^2 P(y_{k-1}/\lambda)}{y_{k-2}} & \text{if $k$ is odd;} \\[.1in] 
\displaystyle\frac{\lambda^2 Q(y_{k-1}/\mu)}{y_{k-2}} & \text{if $k$ is even} 
\end{cases} 
\] 
is a Laurent polynomial in $y_0$ and $y_1$ 
with coefficients in $\AAA=\ZZ[\lambda^{\pm1},\mu^{\pm1},\alpha_i,\beta_i]$. 
This follows from Corollary~\ref{cor:homogeneous-pattern} with $n=2$, 
$P_1=\mu^2 P(x_2/\lambda)$, and $P_2=\lambda^2 Q(x_1/\mu)$. 
} 
\end{example}

\begin{example} 
{\rm 
Consider the sequence $y_0,y_1,\dots$ defined by the recurrence 
\begin{equation} 
y_k=\frac{y_{k-1}^2+c y_{k-1}+d}{y_{k-2}}. 
\end{equation} 
Every term of this sequence is a Laurent polynomial in $y_0$ and $y_1$ 
with coefficients in $\ZZ[c,d]$. 
} 
\end{example} 
 
\begin{example} 
\label{example:rank3-trinomial} 
{\rm 
Define the rational transformations $F_1,F_2,F_3$ by 
% in the variables $x_1,x_2,x_3$ 
\begin{eqnarray} 
\begin{array}{rcrcccl} 
\qquad\qquad 
F_1 : (x_1,x_2,x_3)&\!\!\!\!\mapsto\!\!\!\!&(&\!\!\!\!\displaystyle\frac{x_2+x_3^2+x_2^2 x_3}{x_1},&x_2,&x_3&),\\[.2in] 
F_2 : (x_1,x_2,x_3)&\!\!\!\!\mapsto\!\!\!\!&(&x_1,&\!\!\!\!\displaystyle\frac{x_1+ x_3}{x_2},&x_3&),\\[.2in] 
F_3 : (x_1,x_2,x_3)&\!\!\!\!\mapsto\!\!\!\!&(&x_1,&x_2,&\displaystyle\frac{x_2+x_1^2+x_2^2 x_1}{x_3}\!\!\!\!&). 
\end{array} 
\end{eqnarray} 
Then any composition $F_{i_1}\circ F_{i_2}\circ \cdots$ is given by 
$ 
(x_1,x_2,x_3) \mapsto (G_1,G_2,G_3)$,  
where $G_1,G_2,G_3$ are Laurent polynomials in $x_1,x_2,x_3$ over~$\ZZ$. 
%(This example can be vastly generalized as well.) 
} 
\end{example}

\end{document}